\font\gothic=eufm10.
.
\font\piccolissimo=cmr5.
\font\script=eusm10.
\font\sets=msbm10.
\font\stampatello=cmcsc10.
\font\symbols=msam10.

\def\spaziolungo{\qquad \qquad \qquad \qquad \qquad \qquad }
\def\1{{\bf 1}}
\def\sgn{{\rm sgn}}

\def\avesum{\sum_{x\sim N}}

\def\square{\hbox{\vrule\vbox{\hrule\phantom{s}\hrule}\vrule}}
\def\defineq{\buildrel{def}\over{=}}
\def\defin{\buildrel{def}\over{\Longleftrightarrow}}
\def\doublesum{\mathop{\sum\sum}}

\def\multiplesum{\mathop{\sum \enspace \cdots \enspace \sum}}

\def\C{\hbox{\sets C}}
\def\D{\hbox{\sets D}}
\def\N{\hbox{\sets N}}

\def\R{\hbox{\sets R}}
\def\Z{\hbox{\sets Z}}
\def\Corr{\hbox{\script C}}
\def\EssBdd{\hbox{\symbols n}\,}
\def\modSel{{\widetilde{J}}}

\def\SingSer{\hbox{\gothic S}}
\def\SingInt{{\cal I}}
\def\Res{\mathop{{\rm Res}\,}}
\def\ord{\mathop{{\rm ord}}}
\def\divisor{\hbox{\bf d}}

\par
\centerline{\bf Generations of correlation averages}
\bigskip
\centerline{\stampatello giovanni coppola - maurizio laporta}
\bigskip
{

{\par
{\bf Abstract.} The present paper is a dissertation on the possible consequences
of a conjectural bound for the so-called \thinspace 
modified Selberg integral of the divisor function $d_3$, i.e. a discrete version of the classical Selberg integral, where 
$d_3(n)=\sum_{abc=n}1$ is attached to the Cesaro weight $1-|n-x|/H$ in the short interval $|n-x|\le H$.
Mainly, an immediate consequence is a non-trivial bound for the Selberg integral of $d_3$, 
improving recent results of Ivi\'c based on the standard approach through the moments of the Riemann zeta function on the critical line.  We proceed instead with elementary arguments, by first applying the \lq \lq elementary Dispersion Method\rq \rq\ in order to establish
a link between \lq \lq weighted Selberg integrals\rq \rq \thinspace of any arithmetic function $f$ and averages of correlations of $f$ in short intervals. Moreover, we provide a 
conditional
generalization of our results to the analogous problem on the divisor function $d_k$ for any $k\ge 3$. Further, some remarkable
consequences  on the $2k-$th moments of the Riemann zeta function are discussed. Finally, we also discuss the essential properties that a general function $f$ should satisfy so that the estimation of its Selberg integrals could be approachable by our method. 

\footnote{}{\par \noindent {\it Mathematics Subject Classification} $(2010) : 11{\rm N}37, 11{\rm M}06.$}
}
\bigskip

\par
\noindent {\bf 0. Libretto: introduction and statement of the results.}
\smallskip
\par
\noindent
In the milestone paper [S] Selberg introduced 
a determinant tool in the study of the distribution of prime numbers in {\it short intervals} $[x,x+H]$, i.e. $H=o(x)$ as $x\to \infty$, namely
the integral
$$
\int_{N}^{2N}\Big| \sum_{x<n\le x+H}\Lambda(n)-H\Big|^2 {\rm d}x\ , 
$$
\par
\noindent
where $\Lambda$ is the von Mangoldt function defined as $\Lambda(n)\defineq \log p$ if $n=p^r$ for some prime number $p$ and for some positive integer $r$, otherwise $\Lambda(n)\defineq 0$. Thus, $\Lambda$ 
is a weighted characteristic function 
of the prime numbers and it
is generated by (minus) the logarithmic derivative of the
Riemann zeta function, i.e. its Dirichlet series is 
$-\zeta'(s)/\zeta(s)$ . Further, being a quadratic mean, the Selberg integral  precisely concerns the study of the distribution of primes 
in {\it almost all} short intervals $[x,x+H]$ with at most $o(N)$ exceptional integers $x\in [N,2N]$ as $N\to \infty$. Here we define the {\it Selberg integral}
of any arithmetic function $f:\N \rightarrow \C$ as
$$
J_f(N,H)\defineq \avesum \Big| \sum_{x<n\le x+H}f(n)-M_f(x,H)\Big|^2\, , 
$$
\par
\noindent
where $x\sim N$ means $N<x\le 2N$ and $M_f(x,H)$ is the expected {\it mean value} of $f$ in short intervals (abbreviated as s.i. mean value). 
In order to avoid trivialities, one assumes that the length $H$ of the short interval goes to infinity with $N$.
In view of 
non-trivial bounds of such sums, it is easy to realize that the discrete version $J_{\Lambda}(N,H)$ is
close enough to the original integral introduced by Selberg, so that we feel to be legitimate to use the same symbol for both versions. Such conditions hold for the arithmetic functions we work with and 
the typical case is the $k-$divisor function $d_k$ for $k\ge 3$, where $d_k(n)$ is the number of ways to write $n$ as a product of $k$ positive integer factors (see [C0] and compare $\S3$).
Let us denote the Selberg integral of $d_k$ as
$$
J_k(N,H)\defineq \avesum \Big| \sum_{x<n\le x+H}d_k(n)-M_k(x,H)\Big|^2 
$$
\par
\noindent
with the s.i. mean value of $d_k$ given by 
$$
M_k(x,H)\defineq H\left( P_{k-1}(\log x)+P'_{k-1}(\log x)\right)\, , 
$$
\par
\noindent
where $P_{k-1}$ is the {\it residual polynomial} of degree $k-1$ such that $P_{k-1}(\log x)\defineq \Res_{s=1}\zeta^k(s)x^{s-1}/s$. 
\smallskip
\par				
\noindent
The first author has proved the lower bound $NH\log^4N\ll J_3(N,H)$ for $H\ll N^{1/3-\varepsilon}$ (see [C0]), while, in an attempt to establish a non trivial upper bound, both the authors
have formulated the following conjecture for the so-called {\it modified Selberg integral} of $d_3$,
$$
\modSel_3(N,H)\defineq\avesum \Big| \sum_{0\le |n-x|\le H}\Big( 1-{{|n-x|}\over H}\Big)d_3(n)-M_3(x,H)\Big|^2\ ,
$$
\par
\noindent
where $M_3(x,H)$ is the same s.i. mean value of $J_3(N,H)$. 
\medskip
\par
\noindent {\stampatello Conjecture CL.} {\it If $H\ll N^{1/3}$, then} $\modSel_3(N,H)\EssBdd NH$.
\smallskip
\par
\noindent
Here and in what follows  for convenience we write
$$
A(N,H)\EssBdd B(N,H)
\enspace \hbox{whenever } \enspace
A(N,H)\ll_{\varepsilon} N^{\varepsilon}B(N,H)\quad \forall \varepsilon>0.
$$
\par
\noindent
Moreover, we adopt a further convention on bounds of the {\it width} of $H$, i.e. $\theta \defineq \log H/\log N$ (more in  general 
$\theta$ is defined by $x^{\theta} \ll H\ll x^{\theta}$ for $x\sim N$): the inequality $\theta>\theta_0$ 
(resp. $\theta<\theta_0$) means that there exists a fixed and absolute constant $\delta>0$ such that $\theta \ge \theta_0+\delta$ (resp. $\theta \le \theta_0-\delta$).
In particular, $0<\theta<1$ has to be interpreted as $\delta \le \theta \le 1-\delta$.
\smallskip
\par
As a consequence one has the following result. 
\smallskip
\par
\noindent {\stampatello Theorem 1.} {\it If Conjecture CL holds, then $
J_3(N,H)\EssBdd NH^{3/2}$.}
\medskip
\par
\noindent
Our theorem is an easy deduction by the general link between 
the Selberg integral $J_f(N,H)$ and the corresponding modified one (see \S 4)
$$
\modSel_f(N,H)\defineq \avesum \Big| \sum_{0\le |n-x|\le H}\Big( 1-{{|n-x|}\over H}\Big)f(n)-M_f(x,H)\Big|^2 .
$$
\par
\noindent
Noteworthily Theorem 1 implies an improvement on Ivi\'c's results [Iv2] for $d_3$ both in the bound and in the \lq \lq low\rq \rq \thinspace range where it is valid: while 
Ivi\'c's bound is non-trivial for \enspace $N^{1/6+\delta}\le H\le N^{1-\delta}$, i.e. for width $1/6<\theta<1$, ours is non-trivial for $0<\theta\le 1/3$, that is in the range \enspace $N^{\delta}\le H\le N^{1/3}$. Further, we think that our estimates can be refined in order to get a better range for the width. 
We remark that, still assuming Conjecture CL, the first author [C5] has recently derived the better bound 
$$
J_3(N,H)\EssBdd NH^{6/5}.
$$

\smallskip

Needless to say that our study applies to any divisor function $d_k$, though the conjectured estimates
of the modified Selberg integral,
$$
\modSel_k(N,H)\defineq \avesum \Big| \sum_{0\le |n-x|\le H}\Big( 1-{{|n-x|}\over H}\Big)d_k(n)-M_k(x,H)\Big|^2\ ,
$$
become less and less meaningful as $k$ grows. This is essentially due to the poor state of knowledge about the distribution of $d_k$ in long intervals, namely the known value of the exponent $\alpha_k$ such that (compare $\S3$)
$$
\sum_{n\le x}d_k(n)-xP_{k-1}(\log x)\EssBdd x^{\alpha_k}\ .
\leqno{(\ast)} 
$$
\par
\noindent
We generalize Conjecture CL for $\modSel_k(N,H)$ with $k>3$ as it follows. \smallskip

\noindent {\stampatello General Conjecture CL.} {\it Assume that $\alpha_{k-1}\in [0,1)$ in $(\ast)$ for a fixed integer $k>3$. 
If $H\ll N^{1/2}$, then
$$
\modSel_k(N,H)\EssBdd N^{1-1/k}H^2+NH(N^{1-4/k} + N^{\left(1-1/k\right)\alpha_{k-1}-1/k})\ .
$$
\par				
\noindent
Consequently, setting
$\displaystyle{\widetilde{\theta}_k \defineq \max\Big(1-{4\over k},\Big( 1-{1\over k}\Big)\alpha_{k-1}-{1\over k}\Big)}$\ ,
for every width $\theta\in(\widetilde{\theta}_k,1/2]$
there exists an $\varepsilon_1=\varepsilon_1(\theta,k)>0$ such that}
$$
\modSel_k(N,H)\ll N^{1-4\varepsilon_1}H^2\ .
$$
\par
\noindent
Similarly to the case $k=3$, from the last inequality it follows our second main result.\smallskip
\par
\noindent {\stampatello Theorem 2.} {\it If General Conjecture CL holds, then $
J_k(N,H)\ll N^{1-2\varepsilon_1}H^2$.}
\medskip
\par
\noindent
Again as a consequence one would get an improvement in the low range of $H$ with respect to the results of Ivi\'c for the mean-square of $d_k$ in short intervals [Iv2].
In fact, 
Theorem 1 in [Iv2] holds for $\theta\in(\theta_k,1)$ with $\theta_k$ defined in terms of Carlson's abscissae (see $\S6$). In particular, it holds for
$\theta_4=1/4$, $\theta_5=11/30$, $\theta_6=3/7$, whereas 
we get non-trivial estimates for widths $\widetilde{\theta}_k<\theta\le 1/2$ with $\widetilde{\theta}_4=11/128$, $\widetilde{\theta}_5=1/5$, $\widetilde{\theta}_6=1/3$.
See \S7 for further details, where one finds the so-called $k-${\it folding} method that is
at the core of our conjectures. 
\smallskip
\par
It is apparent from our study that generally $\modSel_f(N,H)$
fits the request of \lq \lq smoothing\rq \rq\ the Selberg integral $J_f(N,H)$, both in the arithmetic and the harmonic analysis 
aspects. The arithmetic matter essentially relies on
a simple observation going back to the Italian mathematician Cesaro around the end of the 19th century:
$$
\sum_{0\le |n-x|\le H}\left( 1-{{|n-x|}\over H}\right)f(n)={1\over H}\sum_{h\le H}\sum_{|n-x|<h}f(n)\ .
$$
\par
\noindent
This is a kind of arithmetic mean of the inner sum in $J_f(N,H)$
and somehow justifies the appearance of the same mean-value term
in the modified Selberg integral. 
The analytic aspects of such a smoothing process will be better understood after the introduction of the 
{\it correlation} $\Corr_f(h)$ in $\S2$, where it is showed that 
Selberg integrals of $f$ are strictly related to averages of $\Corr_f(h)$ in short intervals $|h|\ll H$. 
\par
The next corollaries testify such an intimate link and also conditionally improve recent results [BBMZ] and [IW] on an additive divisor problem for $d_k$. More precisely, they concern the {\it deviation} of $d_k$, i.e.
$$
\D_k(N,H)\defineq\sum_{h\le H}\sum_{n\sim N}d_k(n)d_k(n-h)-{1\over H}\sum_{n\sim N}M_k(n,H)^2=
\sum_{h\le H}\Corr_k(h)-H\sum_{n\sim N}p_{k-1}(\log n)^2
\ ,
$$
\par
\noindent
where $\displaystyle{\Corr_k(h)\defineq\sum_{n\sim N}d_k(n)d_k(n-h)}$ is the correlation of $d_k$ (see \S2) and
$p_{k-1}(\log n)\defineq M_k(n,H)/H$ is
the so-called {\it logarithmic polynomial} of $d_k(n)$.
The following further consequence of Theorem 1 is proved in \S6.
\medskip
\par
\noindent {\stampatello Corollary 1}. {\it Let $N,H$ be positive integers such that $0<\theta=\log H/\log N\le 1/3$. 
Then} 
$$
\D_3(N,H)\EssBdd NH^{3/4} + N^{\alpha_3}H + H^2. 
$$
\medskip
\par
\noindent
Analogously, from Theorem 2 we get a non-trivial estimate for the deviation of $d_k$ for any $k>3$ (see \S7).
\medskip
\par
\noindent {\stampatello Corollary 2.} {\it Under the same hypotheses of Theorem 2 one has, in the same ranges and for the same } $\varepsilon_1$, 
$$
\D_k(N,H)\EssBdd N^{1-\varepsilon_1}H.
$$
\medskip
\par
\noindent
The aforementioned link to [BBMZ] and [IW] results is due to the identity (see $\S3$) 
$$
p_{k-1}(\log x)=P_{k-1}(\log x)+P'_{k-1}(\log x)=\Res_{s=1}\zeta^k(s)x^{s-1}\ , 
$$
\par
\noindent
where one has to be acquainted that our notations $P_{k-1}$ and $p_{k-1}$ are not consistent with those in [BBMZ] and [IW].
In particular, from equation $(3.8)$ in [IW] it turns out that 
$$
H\int_{N}^{2N}p_{k-1}(\log x)^2 {\rm d}x + O_{\varepsilon}\left( N^{1+\varepsilon}\right) 
$$
\par				
\noindent
is the main term in the formul\ae\ for sums of $d_k$ correlations established in [BBMZ] when $k=3$ and in [IW] for every $k\ge 3$. 
Since it is easily seen that
$$
\sum_{n\sim N}p_{k-1}(\log n)^2 - \int_{N}^{2N}p_{k-1}(\log x)^2{\rm d}x\EssBdd 1\ , 
$$
then for every $H\ll N$ it follows 
$$
H\sum_{n\sim N}p_{k-1}(\log n)^2 - H\int_{N}^{2N}p_{k-1}(\log x)^2{\rm d}x\EssBdd N\ , 
$$
\par
\noindent
revealing that within negligible remainders $\EssBdd N$ our $\D_k(N,H)$ is comparable with $\sum_{h\le H}\Delta_k(N;h)$, that is
the average of errors for $d_k$ correlations estimated in [BBMZ] when $k=3$ and in [IW] for every $k\ge 3$. 
Thus, Corollary 1 and the best known $\alpha_3\le 43/96$ (Kolesnik, 1981) imply that $\D_3(N,H) \EssBdd NH^{3/4}$
for $\theta \le 1/3$, which improves [BBMZ] in the low range of short intervals ([IW] bounds are better for $k=3$ when $\theta>1/2$). In fact, their remainders total $\EssBdd N^{13/12}\sqrt H$ , that is worst than $\EssBdd NH^{3/4}$ when $\theta<1/3$, and are non-trivial only for $\theta>1/6$. 
Similarly, since Corollary 2 holds
in the ranges prescribed by the General Conjecture CL, we get an improvement  on the estimation of $\D_k(N,H)$ in the low range of short intervals with respect to [IW] bounds' non-trivial ranges, though they have better high range, say $1/2<\theta<1$.  
\par
The novelty of our approach is that, though conditionally, it leads to valuable improvements with respect to the analogous achievements obtained via the classical  moments of the Riemann zeta function on the critical line $\Re(s)=1/2$, while
we think that the conjectures CL on the modified Selberg integrals might be approachable by
elementary arguments. 
On the other side, estimates of $J_k(N,H)$ have non-trivial consequences on the $2k-$th moments of $\zeta$ (see [C2]) defined as
$$
I_k(T)\defineq \int_{T}^{2T}\Big| \zeta\Big( {1\over 2}+it\Big)\Big|^{2k} {\rm d}t\ .  
$$
Thus, at the moment 
we content ourselves in having found an alternative way to pursue possible improvements on the $2k-$th moments of $\zeta$ at least for relatively low values of $k$. Indeed, in $\S8$ we take a glance at the effect of hypothetical estimates for Selberg integrals on the $2k-$th moments through Theorem 1.1 of [C2], whereas our conjectured $\modSel_3$ bound provides effectively the best known estimate for the $6-$th moment of Riemann $\zeta$ function after recent [C5]. 
In particular, in $\S8$ we prove next result, which
gives a link between conditional bounds of Selberg integrals $J_k$ and bounds of $I_k(T)$.

\medskip

\par
\noindent {\stampatello Theorem 3.} {\it Let  $k\ge 3$ be fixed. If \thinspace 
$\displaystyle{J_k(N,H)\EssBdd N^{1+A}H^{1+B}}$ holds
for $H\ll N^{1-2/k}$ and for some constants $A,B\ge 0$,
then $I_k(T)\EssBdd T^{1+{k\over 2}(A+B)-B}$.}
\medskip
\par
Here, beyond the dependence on $\varepsilon$, the constant involved in $\EssBdd$ may depend on $k$. 
\smallskip
\par
\noindent
As an immediate consequence of Theorem 3 combined with the results in [C5], namely $J_3(N,H)\EssBdd NH^{6/5}$ in the range $H\ll N^{1/3}$, we get 
$$
I_3(T)\EssBdd T^{11/10}\ .
$$
\smallskip
\par
\noindent
This encourages us to follow such a pattern and to pursue non trivial estimates for the modified Selberg integrals of $d_k$ in the future. 

\bigskip

\par
\centerline{\bf Plan of the paper.}

\smallskip

\item{$\S1$} Beyond the aforementioned instances, further variations of the Selberg integral can be considered, according to the weight $w$ suitably attached to an arithmetic function $f$. We give a very short introduction to the so-called $w$-Selberg integral of $f$.
\item{$\S2$} We introduce the correlations of an arithmetic function $f$ and of a weight $w$.
For a wide class of arithmetic functions it is shown through the Dispersion Method that weighted Selberg integrals are strictly related to averages 
of such correlations (see Lemma 1). 
\item{$\S3$} In an attempt to generalize further our results, we have abstracted the essential properties that an arithmetic function $f$ has to satisfy so that
its Selberg integrals can be approachable by our method. Mainly inspired by the prototype $d_k$, we devote this section to 
the definitions and basic properties of the {\it essentially bounded, balanced, quasi-constant} and {\it stable} arithmetic functions.
Inevitably in this analysis one finds references to the famous Selberg Class.
\item{$\S4$} For a real, balanced and essentially bounded function $f$, the second and the third generation of correlation averages in short intervals correspond respectively to the Selberg integral and the modified one.
Here we exploit further the properties of such functions in long intervals to outline a chain of implications that under suitable conditions from a non-trivial estimate for $\modSel_f(N,H)$ generates a non-trivial bound for $J_f(N,H)$. Applying these implications to the divisor functions one gets immediately 
Theorems 1 and 2. Such a bound for $J_f(N,H)$ in turn becomes an effective mean to pursue 
a \lq \lq good deviation\rq \rq\  of $f$, i.e. an error term in the asymptotic formula for the first average of the correlation of $f$ in short intervals,
whereas, as we shortly recall in $\S5$, the expected formula for the single correlation is just conjectural for most significant instances of $f$.  
\item{$\S5$} It is a short excursion on some very special cases 
of single correlations whose conjectured asymptotic formul\ae\ have been proved.
\item{$\S6$} Here one finds the proof of  Corollary 1, that follows rather easily through the general arguments of $\S4$. Corollary 2 follows similarly although it is discussed in the next section.
\item{$\S7$} At least in principle, the strategy applies also to $d_k$ for any $k>3$ as an application of the general $k$-folding method that is described by Lemma 2.
\item{$\S8$} As an application of Theorem 1.1 in [C2], we prove Theorem 3, that emphasizes the consequences of $J_k(N,H)$ bounds on the moments of the Riemann zeta function on the critical line. 
\item{$\S9$} We call upon the Selberg integral of $d_3$ to address
the last word on the best \lq \lq unconditional\rq \rq \thinspace exponent for the $6-$th moment, that becomes an immediate consequence of Conjecture CL,
Theorem 3 and [C5] bounds. 
\medskip
\par
\noindent
\centerline{\bf Some notation and conventions.}
\smallskip
\par
\noindent
If the implicit constants in $O$ and $\ll$ symbols depend on 
some parameters like $\varepsilon>0$, then mostly we specify it by introducing subscripts in such symbols
like $O_\varepsilon$ and $\ll_\varepsilon$, while we avoid subscripts for $\EssBdd$ defined above. 
Notice that the value of $\varepsilon$ may change from statement to
statement, since $\varepsilon>0$ is arbitrarily small.
\par
The relation $f\sim g$ between the functions $f,g$ means that $f=g+o(g)$ as the main variable tends to infinity typically.
No confusion should be possible with the {\it dyadic} notation, 
$x\sim N$, which means that $x$ is an integer of the interval $(N,2N]$, as already said.
\par
The {\it M\"{o}bius function} is defined as $\mu(1)=1$, $\mu(n)= (-1)^r$ if 
$n$ is the product of $r$ distinct primes, and $\mu(n)=0$  otherwise. 
The symbol $\1$ denotes the constant function with value $1$ and $\1_{U}$ is the characteristic function of the set $U$. 
The {\it Dirichlet convolution product} of the arithmetic functions $f_1$ and $f_2$ is 
$$
(f_1 \ast f_2)(n)\defineq \sum_{d|n}f_1(d)f_2(n/d)\quad \forall n\in \N\ . 
$$
\par
\noindent
In particular, we call $f_k\defineq \underbrace{f\ast \cdots \ast f}_{k\, {\rm times}}$ the {\it $k-$fold Dirichlet product} of the arithmetic function $f$.
For any $f,g:\N \rightarrow \C$ the {\it M\"{o}bius inversion formula} states that $f=g\ast \1$ if and only if $g=f\ast \mu$, which is called the {\it Eratosthenes transform} of $f$. For example, 
$\1$ is the Eratosthenes transform of the divisor function $\divisor=\1 \ast \1$. 
More in general, 
$d_k=\underbrace{\1 \ast \cdots \ast \1}_{k\, {\rm times}}=\1_k$ is the $k-$fold Dirichlet product of $\1$ for $k\ge 2$. 
\par
\noindent
For simplicity, in sums like $\sum_{a\le X}$ it is implicit that $a\ge 1$. The distance of a real number $\alpha$ from the nearest integer is denoted by
$\Vert\alpha\Vert=\min(\{ \alpha\}, 1-\{ \alpha\})$, where $\{ \alpha\}$
is the fractional part of $\alpha$. 
As usual, we set $e(\alpha)=e^{2\pi i\alpha}$, $\forall \alpha \in \R$, and $e_q(a)=e(a/q)$, $\forall q\in \N$, $\forall a\in \Z$. 
\medskip
\par
\noindent 

\centerline{\bf Acknowledgements.} 
\smallskip
\par
\noindent 
The authors are very grateful to Alberto Perelli for interesting discussions and invaluable suggestions. 

\vfill
\eject

\par				
\noindent {\bf 1. Preludio: weighted Selberg integrals.}
\smallskip
\par
\noindent
Given positive integers $N$ and $H=o(N)$,
the $w-${\it Selberg integral} of an arithmetic function $f:\N \rightarrow \C$ is the weighted quadratic mean
$$
J_{w,f}(N,H)\defineq \avesum \Big| \sum_n w(n-x)f(n)-M_f(x,w)\Big|^2\ ,
$$
\par
\noindent
where the complex valued {\it weight} $w$ has support in $[-cH,cH]$ for some fixed real number $c>0$, so that
the inner sum is genuinely finite. The term
$M_f(x,w)$ is the expected mean value of $f$ weighted with $w$ in the short interval of length $\ll H$ and its dependence on $w$ has to be specified according to the different instances. In particular, according to the study exposed in $\S3$, when it is possible to define the {\it logarithmic polynomial} $p_f(\log n)$  we set 
$$
M_f(x,w)\defineq \sum_{n}w(n-x)p_f(\log n)=p_f(\log x)\sum_h w(h)+O(H^2 N^{\varepsilon-1})\ . 
$$
\par
\noindent
Clearly, the weighted Selberg integrals include
the most celebrated case of the original Selberg integral, since $J_{\Lambda}(N,H)=J_{u,\Lambda}(N,H)$,
where
$u\defineq \1_{[1,H]}$ is the characteristic function of $[1,H]$. More in general, $J_f(N,H)$ is the
$u-$Selberg integral of $f$.
The modified Selberg integral $\modSel_f(N,H)$, introduced by the first author in [C1], is recognizable as a weighted Selberg integral by taking the {\it Cesaro weight}, say 
$$
C_H(t)\defineq\Big(1-{{|t|}\over H}\Big)_+=\max\Big(1-{{|t|}\over H},0\Big)\ .
$$
\par
\noindent
Since
$$
C_H(t)
={1\over H} \sum_{a\le H-|t|}1 ={1\over H} \doublesum_{{a,b\le H}\atop {b-a=t}}1
={\Corr_u(t)\over H}\ ,
$$ 
\par
\noindent
where $\Corr_u$ is the correlation of the weight $u$ (see next \S2), then we refer to the Cesaro weight as
the \lq \lq normalized correlation\rq \rq \thinspace of $u$. More in general, we {\it smooth} the weighted Selberg integral $J_{w,f}(N,H)$ by defining the {\it modified $w-$Selberg integral} of $f$ as
$$
\modSel_{w,f}(N,H)\defineq J_{\widetilde{w},f}(N,H)\ ,
$$
\par
\noindent
where the new weight $\widetilde{w}$ is the normalized correlation of $w$, i.e.
$$
\widetilde{w}(h)\defineq {1\over H}\doublesum_{{n \thinspace \quad \thinspace m}\atop {n-m=h}}w(n)\overline{w(m)}\ .
$$
\par
\noindent
Another important instance of the weighted Selberg integral has been intensively studied by the first author, i.e. the {\it symmetry integral} of  $f$ given by
$$
J_{\sgn,f}(N,H)\defineq \avesum \Big| \sum_{0\le |n-x|\le H}\sgn(n-x)f(n)\Big|^2\ , 
$$
\par
\noindent
where $\sgn(0)\defineq 0$, $\sgn(t)\defineq |t|/t$ for $t\not=0$, and $M_f(x,\sgn)$ vanishes identically for every $f$.
The study of the symmetry integral has been
motivated by the work 
of Kaczorowski and Perelli (see [KP]), who were 
the very first to exploit a 
strict relation of the classical Selberg integral with
the symmetry properties of the prime numbers. Indeed, in [C], [CS], [C3], [C4], 
the symmetric aspects of the distribution of several samples of arithmetic functions in short intervals
are studied through the analysis of the associated symmetry integral. 
\par
It is worthwhile to exploit the link between $J_{\sgn,f}$ and the {\it modified symmetry integral} $\modSel_{\sgn,f}$ in future papers. As in the case of any odd weight $w$, they demand the s.i. mean values to vanish. On the other side, note that
the normalized correlation of any $w$ is even and inside $M_f\left(x,\widetilde{w}\right)$
we have $\sum_h \widetilde{w}(h)=\left|\sum_a w(a)\right|^2 / H$. 
\par				
Finally, some considerations in $\S3$ make it
plausible that a satisfactory general theory, for the weighted Selberg integrals, may be built within the environment of the Selberg Class. 

\bigskip

\par				
\noindent {\bf 2. Ouverture: weighted Selberg integrals as correlation averages.}
\smallskip
\par
\noindent 
By {\it correlation} of an arithmetic function $f:\N \rightarrow \C$
we mean a shifted convolution sum of the form
$$
\Corr_f(h)\defineq \sum_{n\sim N}f(n)\overline{f(n-h)}
\ , 
$$
\par
\noindent
where the {\it shift} $h$ is an integer such that $|h|\le N$.
Observe that one might take into account only the restriction of $f$ to 
$1\le n\le 3N$. 
Further, a correlation of shift $h$
is essentially a weighted count of the integer solutions $n,m\in(N,2N]$ of the equation $n-m=h$ , as 
$$
\Corr_f(h)
=\doublesum_{{n\sim N \thinspace m\sim N}\atop {n-m=h}}f(n)\overline{f(m)}+O\Big(|h| \cdot \max_{\ell \sim N}|f(\ell)|^2\Big)\ .
$$
\par
\noindent
In the present context it is convenient to define 
the correlation of a weight $w$ by neglecting the $O$-term, i.e. 
$$
\Corr_w(h)\defineq \doublesum_{{a \thinspace \quad \thinspace b}\atop {b-a=h}}w(b)\overline{w(a)}\ . 
$$
\par
\noindent
The reason of such a different definition 
will be clarified after next lemma, where we prove a strict connection between correlations and weighted Selberg integrals
by applying an elementary {\it Dispersion Method}. 
\medskip
\par
\noindent {\stampatello Lemma 1.} {\it Let $N,H$ be positive integers such that 
$H\to \infty$ and $H=o(N)$ as $N\to \infty$. For every 
uniformly bounded weight $w$ with support in $[-cH,cH]$ and every arithmetic function $f$ one has
$$
J_{w,f}(N,H)=\sum_{0\le |h|\le 2cH}\!\!\Corr_w(h)\Corr_f(h)-2\Re\Big(\sum_n f(n)\!\avesum w(n-x)\overline{M_f(x,w)}\Big)+\avesum\left|M_f(x,w)\right|^2
+O\left(H^3\Vert f\Vert_{\infty}^2\right),
$$
\par
\noindent
where $\displaystyle{\Vert f\Vert_{\infty}\defineq \max_{N-cH<n\le 2N+cH}|f(n)|}$\ .
}
\smallskip
\par
\noindent {\stampatello Proof.} By expanding the square and exchanging sums one gets
$$
J_{w,f}(N,H)=\avesum \Big(\sum_n w(n-x)f(n)-M_f(x,w)\Big)\Big(\sum_m \overline{w(m-x)}\cdot \overline{f(m)}-\overline{M_f(x,w)}\Big)
=
$$
$$
= \sum_n f(n)\sum_m \overline{f(m)} \avesum w(n-x)\overline{w(m-x)} 
 -2\Re\Big(\sum_n f(n)\avesum w(n-x)\overline{M_f(x,w)}\Big) 
 +\avesum \left|M_f(x,w)\right|^2\ .
$$
\par
\noindent
Thus, it suffices to show that 
$$
\sum_n f(n)\sum_m \overline{f(m)}\avesum w(n-x)\overline{w(m-x)}
= \sum_h \Corr_w(h)\Corr_f(h)
 + O\left(H^3\Vert f\Vert_{\infty}^2\right)\ , 
$$
\par
\noindent
where we may clearly assume that $a\defineq m-x,b\defineq n-x\in[-cH,cH]$. Consequently, we write
$$
\sum_n f(n)\sum_m \overline{f(m)}\avesum w(n-x)\overline{w(m-x)}=
\sum_{|h|\in[0,2cH]} 
\doublesum_{{n \quad m}\atop {n-m=h}}f(n)\overline{f(m)}
\doublesum_{{a,b\in[-cH,cH]}\atop {{b-a=h}\atop {n-b=m-a\in (N,2N]}}}w(b)\overline{w(a)}\ .
$$
\par
\noindent
Since the condition $n-b=m-a\in (N,2N]$ is implied by $n,m\in (N+cH,2N-cH]$, then the latter is 
$$
\sum_{|h|\in[0,2cH]} 
\doublesum_{{n,m\in (N+cH,2N-cH]}\atop {n-m=h}}f(n)\overline{f(m)}
\doublesum_{{a,b\in[-cH,cH]}\atop {{b-a=h}}}w(b)\overline{w(a)}+
$$
$$				
+O\Big(H\Vert f\Vert_{\infty}^2 \sum_{|h|\in[0,2cH]} 
\Big(
\sum_{n}\sum_{{m\in (N-cH,N+cH]\cup (2N-cH,2N+cH]}\atop {n-m=h}}1+
\sum_{m}\sum_{{n\in (N-cH,N+cH]\cup (2N-cH,2N+cH]}\atop {n-m=h}}1
\Big)
\Big)
= 
$$
$$
=\sum_{|h|\in[0,2cH]} 
\Big(
\doublesum_{{n\sim N,m\in (N+cH,2N-cH]}\atop {n-m=h}}f(n)\overline{f(m)}
 + \doublesum_{{n\in (N+cH,2N-cH],m\sim N}\atop {n-m=h}}f(n)\overline{f(m)}
\Big)
 \doublesum_{{a,b\in[-cH,cH]}\atop {{b-a=h}}}w(b)\overline{w(a)}
 +
$$
$$
+ O\Big( H\Vert f\Vert_{\infty}^2 \sum_{h\ll H}\doublesum_{{m\in (N,N+cH]\cup (2N-cH,2N]}\atop {n-m=\pm h}}1\Big) 
+ O\Big( H\Vert f\Vert_{\infty}^2 \sum_{h\ll H}\doublesum_{{m\in (N-cH,N+cH]\cup (2N-cH,2N+cH]}\atop {n-m=\pm h}}1\Big)
= 
$$
$$
=\sum_{0\le |h|\le 2cH}\doublesum_{{n\sim N \, m\sim N}\atop {n-m=h}}f(n)\overline{f(m)}\doublesum_{{-cH\le a , b\le cH}\atop {b-a=h}}w(b)\overline{w(a)} 
 + O\Big( H\Vert f\Vert_{\infty}^2 \sum_{h\ll H}\sum_{m\in (N-cH,N+cH]\cup (2N-cH,2N+cH]}1\Big)
= 
$$
$$
= \sum_{0\le |h|\le 2cH} \Big(\Corr_f(h)+O(\Vert f\Vert_{\infty}^2 |h|)\Big)\Corr_w(h)
 + O\left(H^2\Vert f\Vert_{\infty}^2 (2cH+1)\right)
= \sum_{0\le |h|\le 2cH} \Corr_w(h)\Corr_f(h)
 + O\left(H^3\Vert f\Vert_{\infty}^2\right).\enspace \square 
$$
\medskip
\par
\noindent {\bf Remark.} The remainder term $O(H^3 \Vert f\Vert_{\infty}^2)$ is essentially due to the estimate of 
\lq \lq short\rq \rq \thinspace segments of length $\ll H$ within \lq \lq long\rq \rq \thinspace sums of length $\gg N$. 
We refer to these short segments as the {\it tails} in the summations. In order to simplify our exposition,
the symbol (T) within some of the following formul\ae\ will warn the reader 
of some tails discarded to abbreviate the formul\ae. 
\par
Thus, by using 
the {\it exponential sum}\footnote{$^1$}{Apart from $\beta$ sign, $\widehat{f}(\beta)$ is also-called the {\it discrete Fourier transform} of $f$. Hereafter we will not specify that it is a finite sum.} 
$$
\widehat{f}(\beta)\defineq \sum_{n\sim N}f(n)e(n\beta)\ ,
$$
\par
\noindent
we write 
$$
\Corr_f(h)=\doublesum_{{\negthinspace m\sim N \thinspace n\sim N}\atop {n-m=h}}f(n)\overline{f(m)}+O(\Vert f\Vert_{\infty}^2 |h|)
=\int_{0}^{1}\big| \widehat{f}(\beta)\big|^2 e(-h\beta){\rm d}\beta + O(\Vert f\Vert_{\infty}^2 |h|)
\buildrel{\hbox{\piccolissimo (T)}}\over{\sim} \int_{0}^{1}\big| \widehat{f}(\beta)\big|^2 e(-h\beta){\rm d}\beta\ . 
$$
\par
\noindent
An important aspect is that the exponential sums, whose coefficients are correlations of a weight $w$,  are non-negative. More precisely, 
$$
\widehat{\Corr_w}(\beta)=\sum_h \Corr_w(h)e(h\beta)=\sum_h \doublesum_{b-a=h}w(b)\overline{w(a)}e(h\beta) =
\Big|\sum_n w(n)e(n\beta)\Big|^2=|\widehat{w}(\beta)|^2\qquad \forall \beta \in [0,1)\ . 
$$
\par
\noindent
In particular, for the correlations of $u=\1_{[1,H]}$ one gets the {\it Fej\'er kernel} 
$$
\widehat{\Corr_u}(\beta)=|\widehat{u}(\beta)|^2=\Big|\sum_n u(n)e(n\beta)\Big|^2
=\Big|\sum_{n\le H}e(n\beta)\Big|^2\ .
$$
\par
\noindent
More in general, the Fej\'er-Riesz Theorem [F] states that any non-negative exponential sum is the square modulus of another exponential sum:
$$
\widehat{w}(\beta)\ge 0\quad \forall \beta \in [0,1) 
\enspace \Longleftrightarrow \enspace \exists v:\
\widehat{w}(\beta)=|\widehat{v}(\beta)|^2\quad \forall \beta \in [0,1)\ . 
$$
\par
\noindent
A particularly easy instance of this theorem follows by recalling that the Cesaro weight is 
the normalized correlation of $u$, i.e.
$C_H(h)=\Corr_u(h)/H$ (see \S1). 
Hence, again Fej\'er's kernel makes its appearance in
$$
\sum_{h} \left(1-{{|h|}\over H}\right)_+e(h\beta)={1\over H}\sum_h \Corr_u(h)e(h\beta)
={|\widehat{u}(\beta)|^2\over H}\ .
$$
\par				
\noindent
We also use to say that the Cesaro weights are {\it positive definite}. We think that
basically  such a property makes the aforementioned smoothing process work for the modified Selberg integral under suitable hypotheses on the function $f$ (see $\S3$), 
while for
the classical Selberg integral $J_{f}(N,H)$ and the symmetry integral $J_{\sgn,f}(N,H)$ it is plain that
\lq \lq $u$\rq \rq \thinspace and \lq \lq $\sgn$\rq \rq \thinspace are far from being positive definite weights. 
One could exploit such a positivity condition in order to prove a non-trivial result for the 
modified Selberg integral of the divisor function $d_k$. 
As a general strategy, from a non positive definite weight $w$ 
with support of length at most $H$ one could call for
its normalized correlation $\displaystyle{\widetilde{w}={\Corr_w/H}}$ generating the non negative exponential sum
$\displaystyle{\widehat{\Corr_w}(\beta)/H=|\widehat{w}(\beta)|^2/H}$. 

\bigskip

\par
\noindent {\bf 3. Starring: essentially bounded, balanced, quasi-constant and stable arithmetic functions.}
\smallskip
\par
\noindent
The wide class of arithmetic functions under our consideration 
consists of functions bounded asymptotically by every arbitrarily 
small power of the variable according to the definition\footnote{$^2$}
{That is $f$ satisfies one of the {\it Selberg class} axioms, the so-called {\it Ramanujan hypothesis} (see esp. [De]).}:
$$
f\enspace \hbox{\rm is} \enspace \hbox{\stampatello essentially \thinspace bounded} 
\enspace \defin \enspace 
\forall \varepsilon>0 \quad f(n)\ll_{\varepsilon} n^{\varepsilon}\ ,
$$
\par
\noindent
that we shortly denote by writing $f\EssBdd 1$. However,
in several circumstances one has to deal with arithmetic
functions having support in an interval of length $\ll N$ or simply with functions restricted to such an interval. Thus, more in general among the essentially bounded functions
 we include $f$ such that $f(n)\ll_{\varepsilon} N^{\varepsilon}\ \forall \varepsilon>0$.
\par
\noindent
A well-known prototype of an essentially bounded function is the divisor function $d_k$, whose Dirichlet series is
$\zeta(s)^k$. 
Similarly, the Dirichlet series 
$$
F(s)\defineq \sum_{n=1}^{\infty}{{f(n)}\over{n^s}}
$$
\par
\noindent
is defined at least in the right half-plane $\Re(s)>1$, whenever the generating function $f$ is essentially bounded
(say the abscissa of absolute convergence is $\sigma_{ac}\le 1$). 
Recall that through Perron's formula the 
expansion of $F$ at $s=1$ leads to an asymptotic formula for the summation function from   
$$
\sum_{n\le x}f(n)={1\over {2\pi i}}\int_{c-i\infty}^{c+i\infty}F(s){{x^s}\over s}{\rm d}s\quad (x\not \in \N)\ , 
$$
\par
\noindent
where $c>\max(0,\sigma_{ac})$. More precisely, with the aid of further properties 
of the Dirichlet series and the Residues Theorem, such an asymptotic formula becomes
$$
\sum_{n\le x}f(n)=\Res_{s=1}F(s){{x^s}\over s}+R_f(x)\ ,
$$
\par
\noindent
where $R_f(x)$ is an error term as long as it is smaller than the main term. 
Assuming that $F$ is meromorphic and denoting the {\it polar order}\footnote{$^3$}{That is the order of the pole of $F$ at $s=1$.} of $F$ by 
$m_F\defineq\ord_{s=1}F$, the main term is more explicitly written as (compare [De]) 
$$
x P_f(\log x)\defineq \Res_{s=1}F(s){{x^s}\over s}\ , 
$$
\par
\noindent
where $P_f$ is the {\it residual polynomial} of $f$, which has degree $m_F-1$, while
it vanishes identically when $m_F< 1$. 
For the remainder term a good estimate would be (compare [De] again) 
$$
R_f(x)\EssBdd x^{\alpha(f)}
$$
\par
\noindent
with a suitable $0\le \alpha(f)<1$ (negative values of exponent $\alpha(f)$ are possible, but discarded as \lq \lq meaningless\rq \rq).
\par
\noindent
This is the case for any divisor function $d_k$. Indeed, from $(\ast)$ of \S0 one has
$$
\sum_{n\le x}d_k(n)=xP_{k-1}(\log x)+\Delta_k(x)
\enspace \enspace \hbox{\rm with} \enspace \enspace 
\Delta_k(x)\EssBdd x^{\alpha_k}, 
$$
\par				
\noindent
where the degree of the residual polynomial $P_{k-1}$ (see \S0) is $k-1$, because the polar order of $\zeta^k$ is
$m_k=k$,  and 
$\alpha_k\le 1-1/k$ is what one can infer inductively from the elementary Dirichlet hyperbola method applied to the first
case $k=2$. More precisely, one has $\Delta_k(x)\ll x^{1-1/k}\log^{k-2}x$. 
Now, by partial summation it is easy to determine the {\it logarithmic polynomial} $p_{k-1}$ such that 
$$
\sum_{n\le x}p_{k-1}(\log n)=x P_{k-1}(\log x)+O(\log^{k-1}x)\ .
$$
\par
\noindent
Thus, we get the decomposition 
$$
d_k(n)\defineq p_{k-1}(\log n)+\widetilde{d_k}(n)\ ,
$$
\par
\noindent
that is a \lq \lq balancing\rq \rq \ of $d_k(n)$ from which the very slowly increasing 
polynomial $p_{k-1}(\log n)$ is subtracted. 
The Dirichlet series generated by $\widetilde{d_k}$ times $x^s/s$ has zero residue at $s=1$. 
Moreover, $(\ast)$ is equivalent to
$$
\sum_{n\le x}\widetilde{d_k}(n)
\EssBdd x^{\alpha_k}.
\leqno(\widetilde{\ast})
$$
\par
\noindent
This invites to formulate the following definitions\footnote{$^4$}{Although with a different meaning, such a
terminology has been coined by Ben Green and Terence Tao. Mainly, Green [Gr] calls {\it balanced} a function $f-\delta$ when $f$ is a characteristic function of a set with density $\delta$.}: 
$$
f \enspace \hbox{is \stampatello balanced} 
\enspace \defin \enspace \Res_{s=1}{{x^s}\over s}F(s)=
\Res_{s=1}{{x^s}\over s}\sum_{n=1}^{\infty}{{f(n)}\over{n^s}}=0, \enspace \forall x 
$$
\par
\noindent
(that is $F$ has an analytic continuation in $s=1$, because $m_F\le 0$), 
$$
f \enspace \hbox{is \stampatello well-balanced of\thinspace exponent } \alpha 
\enspace \defin \enspace 
\sum_{n\le x}f(n)\EssBdd x^{\alpha}\ \hbox{for some}\ \alpha\in[0,1)\ . 
$$
\par
\noindent
Of course, any well-balanced function is also balanced because the previous bound implies that the Dirichlet series is regular at $s=1$. 
However, the converse needs not to be true, as $\Lambda(n)-1$ is a balanced function, but the existence of an exponent $\alpha<1$ is still far from being proved (see some further comments below). 
\smallskip
\par
An essentially bounded
 arithmetic function $a:\N \rightarrow \C$ is said to be {\stampatello quasi-constant}
if there exists  $A\in C^1([1,+\infty),\C)$ such that\footnote{$^5$}{The condition on the derivative $A'$ implies that $A$ is essentially bounded, provided $a$ depends only on $n$. However, we leave open the possibility that $a$ and $A$ might depend on auxiliary parameters. } $A'(t)\EssBdd 1/t$
and $\left.A\right|_{\N}=a$.
Clearly, the logarithmic polynomial $p_{k-1}(\log n)$ 
is quasi-constant with respect to $n$ and this, together with the fact that $\widetilde{d_k}$ is a well-balanced 
arithmetic function of exponent $\alpha_k$, suggests the following further definition.
\smallskip
\par
An arithmetic function  $f$ is {\stampatello stable of exponent} $\alpha$ 
if there exist a quasi-constant function $a$ and a well-balanced function
$b$ of exponent $\alpha$ such that $f=a+b$,
while the {\stampatello amplitude} of $f$ is defined as
$$
\alpha(f)\defineq \inf \, \{ \alpha \in (0,1) : f\ \hbox{is  stable of exponent } \alpha\}\ .
$$
\par				
\noindent
Recall that the Dirichlet divisor problem requires to prove the conjectured amplitude $\alpha_2=\alpha(\divisor)=1/4$,
while one infers  $\alpha_2\le 1/2$ by the Dirichlet hyperbola method and 
the best known bound at the moment is 
$\alpha_2\le 141/416$ (Huxley, 2003). In what follows, $\alpha_k=\alpha(d_k)$ is the best possible exponent in $(\ast)$ and $(\widetilde{\ast})$. 
\medskip
\par
According to Ivi\'c [Iv2], the mean value in the Selberg integral of any arithmetic function $f$, whose Dirichlet series $F(s)$ converges absolutely at least in the half-plane $\Re(s)>1$ and is meromorphic in $\C$, has the {\it analytic form} given by
$$
M_f(x,H)\defineq H\left( P_f(\log x)+P'_f(\log x)\right), 
$$
\par				
\noindent
where $P'_f$ is the derivative of the residual polynomial of $f$. 
We remark that $M_f(x,H)$ is linear in $f$ and
is {\it separable}, i.e. the variables $H, x$ are separated. 
Recall that $p_f\defineq P_f+P'_f$ is the logarithmic polynomial. 
\par
Philosophically speaking, although completely justified from an analytic point of view (using mean-value theorem and derivatives bounds), such a choice of $M_f(x,H)$  is not satisfactory, since one should expect to find the same Selberg integral for $f$ and its balanced part $\widetilde{f}=f-p_f$. 
Indeed, this is the case whenever we define the s.i. mean-value as (see $\S1$) 
$$
\sum_{x<n\le x+H}p_f(\log n)=Hp_f(\log x)+O_{\varepsilon}(N^{\varepsilon}H^2/N). 
$$
\par
\noindent
Of course, up to negligible remainders this is still possible under Ivi\'c's hypothesis. Let us give an idea of a possible extension of these considerations to the case of a more general function $f$. 
\medskip
\par
\noindent
Bearing in mind $(\ast)$, given any arithmetic function $f$, we call a polynomial $P_f$ such that
$$
\sum_{n\le x}f(n)=xP_f(\log x)+O_{\varepsilon}\left( x^{\varepsilon+\alpha}\right)\quad \hbox{with}\ \alpha<1
\leqno{(\ast)_f}
$$
\par
\noindent
the {\it residual polynomial} of $f$, although
$P_f$ is not necessarily defined from
the residues in $s=1$ with the Dirichlet series $F(s)$. 
Then, let us define the {\it logarithmic polynomial} of $f$ as
$$
p_f(\log x)\defineq{d\over {dx}}\left( xP_f(\log x)\right)=P_f(\log x)+P'_f(\log x)\ .
$$
\par
\noindent
Under Ivi\'c's hypothesis on $f$ this immediately yields
$$
M_f(x,H)=Hp_f(\log x) 
$$
\par
\noindent
and
$$
p_f(\log x)=\Res_{s=1}F(s){{x^{s-1}}\over s}+\Res_{s=1}F(s){{(s-1)x^{s-1}}\over s}=\Res_{s=1}F(s)x^{s-1}\ .
$$
\par
\noindent
In the case $f=d_k$ this property allows us to
to compare [IW] results to ours (see \S0). 
\par
On the other side, by assuming the sole property $(\ast)_f$ one gets a unique polynomial $p_f$ such that 
$$
\sum_{n\le x}p_f(\log n)=\int_{1}^{x}p_f(\log t){\rm d}t+O(\log^{c}x)=x P_f(\log x)+O(\log^{c}x)\sim
\sum_{n\le x}f(n)\ ,
$$
\par
\noindent
where $c\ge 0$ is the degree of $p_f$. Since $p_f(\log x)$ is a quasi-constant function,
this implies that
$$
\sum_{x<n\le x+H}f(n)\sim \sum_{x<n\le x+H}p_f(\log n)\sim Hp_f(\log x)\ .
$$
\par
\noindent
Thus every arithmetic function $f$ satisfying $(\ast)_f$ admits a logarithmic polynomial $p_f$ and
the {\it analytic form} of the mean-value in short intervals, inside the $w$-Selberg integral, of such a function $f$ is 
$$
\sum_{n}w(n-x)p_f(\log n)=p_f(\log x)\sum_h w(h)+O(H^2 L^{c-1}/N), 
$$
where $L\defineq \log N$ hereafter.
In particular, the analytic form of the mean value in the Selberg integral of $d_3$ is explicitly given by  
$$
M_3(x,h)=h\big(P_2(\log x)+P'_2(\log x)\big)=h\Big( {1\over 2}\log^2 x + 3\gamma \log x + 3\gamma^2+3\gamma_1\Big)\ ,
$$
\par				
\noindent
with 
$\displaystyle{
P_2(t)=t^2/2+(3\gamma-1)t+(3\gamma^2+3\gamma_1-3\gamma+1)}$, where $\gamma$ is the Euler-Mascheroni constant and $\gamma_1$ is a Stieltjes constant defined as
$$
\gamma_1 \defineq \lim_{m}\Big({{\log^2 m}\over 2}-\sum_{j\le m}{{\log j}\over j}\Big)\ .
$$
\par
\noindent
Recall that this is also related to the summation formula [Ti] 
$$
\sum_{n\le x}d_3(n)=xP_2(\log x)+O(x^{2/3}\log x)\ . 
$$
\par
\noindent
From the application of the $3-${\it folding} method (see $\S7$), recalling that $x\sim N$,
it comes out  that the {\it arithmetic form} of the mean value in the Selberg integral of $d_3$ is 
(here $M=[(N-h)^{1/3}]$, compare \S7) 
$$
\widetilde{M}_3(x,h)\defineq h\Big(
\sum_{q\le {x\over M}}{{\divisor(q)}\over q} 
 + \sum_{d_1<M}{1\over {d_1}}\sum_{d_2\le {x\over {d_1 M}}}{1\over {d_2}} 
 + \Big(\sum_{d<M}{1\over d}\Big)^2\Big)\ . 
$$
\par
\noindent
Indeed, let us show that $M_3(x,h)$ can be replaced by $\widetilde{M}_3(x,h)$ within
$\modSel_3(N,H)$ at the cost of a negligible error term. More precisely, we prove that
$$
\widetilde{M}_3(x,h)-M_3(x,h)
\EssBdd hN^{-1/3}\qquad \hbox{uniformly}\ \forall x\sim N\ .
$$
\par
\noindent
At this aim, we apply Amitsur's formula [A] with Tull's error term [Tu]\hskip0cm
\footnote{$^6$}{ \hskip-0.3cm
Amitsur derived a symbolic method to calculate main terms of asymptotic formul\ae.  Tull, a student of Bateman, gave a refined partial summation, that allows here to transfer error terms from the formula for $\sum_{q\le Q}\divisor(q)$,  like Dirichlet's classical $O(\sqrt Q)$, to this formula for $\sum_{q\le Q}\divisor(q)/q$.}, i.e.
$$
\sum_{q\le Q}{{\divisor(q)}\over q}={{\log^2 Q}\over 2}+2\gamma \log Q + (\gamma^2+2\gamma_1) + O\Big({1\over {\sqrt Q}}\Big)\ ,
$$
\par
\noindent
to get 
$$
\sum_{q\le {x\over M}}{{\divisor(q)}\over q}
={1\over 2}\log^2 x 
- (\log M)\log x + {1\over 2}\log^2 M + 2\gamma \log x - 2\gamma \log M + (\gamma^2 + 2\gamma_1) 
+ O(N^{-1/3})\ . 
\leqno{i)}
$$
\par
\noindent
From the standard asymptotic formula for $\displaystyle{\sum_{d<D}{1\over d}}$ one has
$$
\sum_{d_1<M}{1\over {d_1}}\sum_{d_2\le {x\over {d_1 M}}}{1\over {d_2}}
=\left(\log M + \gamma\right)\log x + \gamma^2 - \log^2 M - {{\log^2 M}\over 2} + \gamma_1 
+ O(N^{-1/3}L)\ , 
\leqno{ii)}
$$
$$				
\Big(\sum_{d<M}{1\over {d}}\Big)^2
=\log^2 M+2\gamma \log M + \gamma^2 
+ O(N^{-1/3}L)\ . 
\leqno{iii)}
$$
\par
\noindent
Thus, $i)$, $ii)$ and $iii)$ imply the claimed inequality, because
$$
{{\widetilde{M}_3(x,h)}\over h} = {1\over 2}\log^2 x 
 + 3\gamma \log x 
 + (3\gamma^2+3\gamma_1) 
 + O(N^{-1/3}L)
= {{M_3(x,h)}\over h} + O_{\varepsilon}(N^{\varepsilon-1/3})\ . 
$$
\par
\noindent
Such a proximity of the two terms $\widetilde{M}_3(x,h)$ and $M_3(x,h)$ suggests that, 
even for more general essentially bounded function $f$, one should expect the same mean value term $M_f(x,h)$
in short intervals for the Selberg integral and for the modified one, whenever
a suitable arithmetic form
$M_f(x,h)\approx h\sum_q g(q)/q$, with $f=g\ast\1$,
is proved to be sufficiently close to the analytic form (determined by the residues 
of the Dirichlet series generated by $f$). This seems to be reliable at least when $M_f(x,t)$ is
{\it separable}, i.e. $M_f(x,t)=t{\cal M}_f(x)$, 
with $t=o(x)$and  ${\cal M}_f(x)$ is a slowly varying function with respect to $x$ (namely, 
a small ${\cal M}'_f(x)$ like the $x-$derivative of polynomials in the variable $\log x$). Indeed, the identities 
$$
\sum_{0\le |n-x|\le h}\Big( 1-{{|n-x|}\over h}\Big)f(n)={1\over h}\sum_{m\le h}\sum_{0\le |n-x|<m}f(n)
\enspace , \enspace 
\sum_{0\le |n-x|<m}f(n) \approx M_f(x,2m-1)= (2m-1){\cal M}_f(x)
$$
\par
\noindent
imply together 
$$
\sum_{0\le |n-x|\le h}\Big( 1-{{|n-x|}\over h}\Big)f(n) \approx {1\over h}\sum_{m\le h}M_f(x,2m-1) 
= h {\cal M}_f(x) = M_f(x,h)\ .
$$
\par
\noindent
We refer the reader to the further discussion in $\S 7$
for the generalization to any divisor function $d_k$ through the so-called $k-${\it folding} method. 
Of course, $d_k$ is not the only function suitable for $(\ast)_f$. 
For example, De Roton [De] has showed that, if the Dirichlet series $F$ belongs to the so-called {\it Extended Selberg Class} (ESC)
with\footnote{$^7$}{The degree of $F$ in ESC is defined in terms of its functional equation (see [De] for details). } ${\rm deg}\, F\ge 1$, then
$(\ast)_f$ holds with $\displaystyle{\alpha={{{\rm deg}\, F-1}\over {{\rm deg}\, F+1}}}$.
Further, if $f$ is also multiplicative (so that its Dirichlet series $F$ has a suitable Euler product) and essentially bounded, then $F$ belongs to
the special subset of  ESC, called {\it Selberg Class} (see [KP(012)]). Hence, according to our definitions De Roton's result (after [KP(012)] breakthrough on Selberg Class) becomes:
\smallskip
\par
\centerline{$f$ has Dirichlet series in the Selberg Class with degree $d\ge 1$
\enspace $\Longrightarrow$ \enspace 
$f$ is stable of exponent $\displaystyle{{d-1}\over {d+1}}$\ .}
\smallskip
\par
\noindent
In particular, this applies to
any $d_k$, since every power $\zeta^k$ for $k\ge 1$ belongs to the Selberg Class with ${\rm deg}\,\zeta^k=k$. 
However, the De Roton exponent $\displaystyle{{k-1}\over {k+1}}$ is weaker than the one obtained by other methods. 
\medskip
\par
\noindent
Actually, the bound $\alpha_k=\alpha(d_k)\le 1/2$ for $k\le 4$ assures that the function $d_k$ for $k\le 4$ admits (at least) square-root cancellation for the error terms, a property shared by every stable arithmetic function with a sufficiently small amplitude. This motivates the following further definition: 
$$
f\ \hbox{is \stampatello random} 
\enspace \defin \enspace 
f\ \hbox{\rm is\enspace stable\enspace of\enspace amplitude}\ \alpha(f)\le 1/2\ . 
$$
\par
\noindent
For example, it is well-known (see [D]) that the Riemann Hypothesis (RH) is
equivalent to  the inequality
$$
R_{\Lambda}(x)\defineq\sum_{n\le x}(\Lambda(n)-1)\EssBdd \sqrt x\ ,
$$
\par
\noindent
that is to say the von Mangoldt function is random\footnote{$^8$}{
It is also a well-known fact the equivalence between RH and the {\it randomness}
of the M\"obius function
$\mu$, suggesting that $\mu$ behaves like $\Lambda-1$.
We refer to [IwKo] for further details on
the {\it M\"obius randomness law}.}.
On the other side, it is also well-known that unconditionally the inequality 
$R_{\Lambda}(x)\EssBdd x^\alpha$ holds only if  
$\alpha=\alpha(\Lambda)\ge 1/2$. In other words, $\Lambda$ cannot be stable of exponent $\alpha<1/2$ . Such a circumstance  is better expressed by the definition:
$$
f \enspace \hbox{is \stampatello strictly random} 
\enspace \defin \enspace 
f\ \hbox{\rm is\enspace stable\enspace of\enspace amplitude}\ \alpha(f)=1/2\ . 
$$
\par
\noindent
Thus, RH is equivalent to say that prime numbers are {\it strictly random}. The existence of $\alpha(\Lambda)<1$ corresponds to a quasi-RH because
of another well-known analytic property of the prime numbers\footnote{$^9$}{Recall also that $\alpha(\mu)=\alpha(\Lambda)$.}:
$$
\alpha(\Lambda)=\sup\{\beta : \zeta(\beta+i\gamma)=0\ \hbox{\rm for\enspace some}\ \gamma \neq 0\}\ .
$$
\par				
The Polya-Vinogradov inequality (see [D]) provides 
a non-conditional example of a strictly random arithmetic function, namely any non-principal Dirichlet character $\chi(\bmod \; q)$, since it yields 
\vskip-0.1cm
$$
\sum_{n\le q}\chi(n)\ll \sqrt{q}\log q
\; 
\EssBdd q^{1/2}\ ,
$$
\vskip-0.1cm
\par
\noindent
which is known to be essentially optimal (compare [Go] \& [Te]). The actual results and the expected values $\alpha_k = (k-1)/(2k)$ in general (see [Iv0]) reveal that the $k-$divisor functions are not strictly random.
\par
Returning back to our integrals, if $f$ is real and essentially bounded, then from Lemma 1 we have
\vskip-0.1cm
$$ 
J_f(N,H)\buildrel{\hbox{\piccolissimo (T)}}\over{\sim}\sum_h \Corr_u(h)\Corr_f(h)-2\sum_n f(n)\avesum u(n-x)M_f(x,H)+\avesum M_f^2(x,H)\ , 
$$
\vskip-0.1cm
$$
\modSel_f(N,H)\buildrel{\hbox{\piccolissimo (T)}}\over{\sim}\sum_h \Corr_{\Corr_u/H}(h)\Corr_f(h)-2\sum_n f(n)\avesum {{\Corr_u(n-x)}\over H}M_f(x,H)+\avesum M_f^2(x,H)\ . 
$$
\vskip-0.1cm
\par
\noindent
When $f$ is also balanced, then $M_f(x,H)$ vanishes identically. Consequently, 
$$ 
J_f(N,H)\buildrel{\hbox{\piccolissimo (T)}}\over{\sim}\sum_h \Corr_u(h)\Corr_f(h)
\qquad \hbox{and} \qquad 
\modSel_f(N,H)\buildrel{\hbox{\piccolissimo (T)}}\over{\sim}\sum_h \Corr_{\Corr_u/H}(h)\Corr_f(h)\ . 
$$

\bigskip

\par
\noindent {\bf 4. Story: smoothing correlations by arithmetic means.}
\smallskip
\par
\noindent
Recalling that 
$$
\Corr_f(h)
\buildrel{\hbox{\piccolissimo (T)}}\over{\sim} \int_{0}^{1}\big| \widehat{f}(\beta)\big|^2 e(-h\beta){\rm d}\beta\qquad 
\hbox{with }\enspace \displaystyle{\widehat{f}(\beta)=\sum_{N<n\le 2N}f(n)e(n\beta)}\ ,
$$
\par
\noindent
one easily infers
$$
\sum_{h}u(h)\Corr_f(h)=\int_{0}^{1}\big| \widehat{f}(\beta)\big|^2 \widehat{u}(-\beta){\rm d}\beta + O(H^2 \Vert f\Vert_{\infty}^2)
\buildrel{\hbox{\piccolissimo (T)}}\over{\sim} \int_{0}^{1}\big| \widehat{f}(\beta)\big|^2 \widehat{u}(-\beta){\rm d}\beta 
\leqno{(I)}
$$
\vskip-0.5cm
$$
\sum_{h}\Corr_u(h)\Corr_f(h)=\int_{0}^{1}\big| \widehat{f}(\beta)\big|^2 \widehat{\Corr_u}(-\beta){\rm d}\beta + O(H^3 \Vert f\Vert_{\infty}^2)
\buildrel{\hbox{\piccolissimo (T)}}\over{\sim} \int_{0}^{1}\big| \widehat{f}(\beta)\big|^2 \big|\widehat{u}(\beta)\big|^2 {\rm d}\beta 
\leqno{(II)}
$$
\vskip-0.5cm
$$
\sum_{h}\Corr_{\Corr_u/H}(h)\Corr_f(h)=\int_{0}^{1}\big| \widehat{f}(\beta)\big|^2 {{\big|\widehat{\Corr_u}(\beta)\big|^2}\over {H^2}}{\rm d}\beta + O(H^3 \Vert f\Vert_{\infty}^2)
\buildrel{\hbox{\piccolissimo (T)}}\over{\sim} \int_{0}^{1}\big| \widehat{f}(\beta)\big|^2 \big|\widehat{u}(\beta)\big|^2 \cdot {{\big|\widehat{u}(\beta)\big|^2}\over {H^2}}{\rm d}\beta 
\leqno{(III)}
$$
\par
\noindent
In particular, for every balanced real function $f\EssBdd 1$, from the previous section we get
$$
J_f(N,H)\buildrel{\hbox{\piccolissimo (T)}}\over{\sim} \int_{0}^{1}\big| \widehat{f}(\beta)\big|^2 \big|\widehat{u}(\beta)\big|^2 {\rm d}\beta 
\qquad \hbox{and} \qquad 
\modSel_f(N,H)\buildrel{\hbox{\piccolissimo (T)}}\over{\sim} \int_{0}^{1}\big| \widehat{f}(\beta)\big|^2 \cdot {{\big|\widehat{u}(\beta)\big|^4}\over {H^2}}{\rm d}\beta\ . 
$$
\par
\noindent
Formulae (I), (II) and (III) correspond respectively to the following iterations of correlations' averages.
$$
\sum_{h\le H}\Corr_f(h)=\sum_{h}u(h)\Corr_f(h)\qquad {\hbox{(sums of correlations)}}
\leqno{1\hbox{\stampatello st generation}:} 
$$
$$
\sum_{h\le H}\sum_{0\le |a|<h}\Corr_f(a)
 =\sum_{h}\Corr_u(h)\Corr_f(h)\qquad \hbox{(double sums)}
\leqno{2\hbox{\stampatello nd generation:}}
$$
$$
\sum_{h}\doublesum_{h_2-h_1=h}{{\Corr_u(h_1)\Corr_u(h_2)}\over H^2}\Corr_f(h)
=\sum_{h}\Corr_{\Corr_u/H}(h)\Corr_f(h)\qquad  \hbox{(average of double sums)}
\leqno{3\hbox{\stampatello rd generation:}}
$$
\par				
\noindent
Such an obstinate process of averaging is motivated by the fact that
it is rarely possible to find an asymptotic formula for the single correlation $\Corr_f$ when $f$ is
a significant arithmetic function. As we already said, the correlation of $f$ counts the number of $h-$twins 
not only when $f$ is a pure characteristic function (the von Mangoldt 
function is a typical case). In general, the underlying Diophantine equation is a binary problem that is out of reach with the present methods
(see next $\S5$ for some details). On the other side, the higher is the degree of a generation of the correlations' averages, the smoother are such averages and consequently we have more hope to get non-trivial asymptotic estimates. However, even at a $2$nd generation level this hope is quite frustrated by the lack
of efficient elementary methods\footnote{$^{10}$}
{This will be coped by our forthcoming paper on mean-squares in short intervals ($w-$Selberg integrals).} 
 to bound directly the Selberg integral. Indeed, Ivi\'c [Iv2] applies Riemann zeta moments
since the Selberg integral of $d_k$ has a strong connection with them (see $\S8$).
Further, it is interesting to analyze the cost of the loss when a
non-trivial information on the correlations' averages at some $n$th
generation level is transfered to the averages of the $(n-1)$th 
generation. For example, assuming that $f$ is real, essentially bounded and balanced,  
from the trivial inequality $\widehat{u}\ll H$ one immediately has 
$$
\modSel_f(N,H)\EssBdd J_f(N,H) + H^3\ . 
$$
\par
\noindent
Then, in order to obtain an inequality in the opposite direction,
we appeal to the formul\ae\ deduced from (II) and (III) and write
$$
J_f(N,H)\EssBdd\!\int_{0}^{1}\!\big|\widehat{f}(\beta)\big|^2 \big|\widehat{u}(\beta)\big|^2 {\rm d}\beta + H^3 
\EssBdd \sqrt{N\!\!\int_{0}^{1}\!\big| \widehat{f}(\beta)\big|^2\big|\widehat{u}(\beta)\big|^4 {\rm d}\beta} + H^3\EssBdd 
\sqrt{NH^2 \modSel_f(N,H)}+\sqrt{N}H^{5/2}+H^3,
$$
\par
\noindent
where we have applied the Cauchy inequality and the Parseval identity (with $f\EssBdd 1$). Thus, if
for some $H$ one has a non-trivial estimate of the kind 
$$
\modSel_f(N,H)\ll NH^2/G
$$
\par
\noindent
with some gain $G\to \infty$, then the previous formula implies (for the same range of $H$)
$$
J_f(N,H)\EssBdd NH^2G^{-1/2}+N^{1/2}H^{5/2}\ .  
$$
\par
\noindent
This gives Theorem 1 by taking $G=H$ in Conjecture CL.
\par
Hence, by the sole application of the Cauchy inequality a third generation gain $G\ll N/H$ leads to the gain $\sqrt G$ 
for the second generation estimate. We say that the {\it exponent's gain has halved}.
\par
The same phenomenon occurs for a general weight $w$ with the
alternative approach that we describe here assuming $\modSel_{w,f}(N,H)\ll NH^2/G$. 
By taking $E\defineq \{\beta:\ \big|\widehat{w}(\beta)\big|\ge\varepsilon H\}$ we have that
the following 
$$
\int_{0}^{1}\big| \widehat{f}(\beta)\big|^2 \big|\widehat{w}(\beta)\big|^2 {\rm d}\beta\ll \varepsilon^2H^2 
\int_{[0,1)\setminus E}\big| \widehat{f}(\beta)\big|^2 {\rm d}\beta+{1\over \varepsilon^2}
\int_{E}\big| \widehat{f}(\beta)\big|^2 \cdot {{\big|\widehat{w}(\beta)\big|^4}\over {H^2}}{\rm d}\beta
$$
$$
\ll NH^2||f||_\infty^2\varepsilon^2+{1\over \varepsilon^2}\widetilde{J}_{w,f}(N,H)\ll NH^2||f||_\infty^2(\varepsilon^2+G^{-1}\varepsilon^{-2})
$$
\par
\noindent
is an optimal bound when $\varepsilon=G^{-1/4}$ and we get the same halving of the exponent's gain as before.
\par
Now, recallling that  $|\widehat{u}(\beta)|\gg H$ when $|\beta|<1/(2H)$ (see [D], Ch.25), we get also  
$$
H^2\int_{-{1\over {2H}}}^{{1\over {2H}}}\Big| \widehat{f}(\beta)\Big|^2 {\rm d}\beta
\ll \int_{-{1\over {2H}}}^{{1\over {2H}}}\Big| \widehat{f}(\beta)\Big|^2\cdot {{\big|\widehat{u}(\beta)\big|^4}\over {H^2}} {\rm d}\beta
\le
\int_{0}^{1}\big| \widehat{f}(\beta)\big|^2 \cdot {{\big|\widehat{u}(\beta)\big|^4}\over {H^2}}{\rm d}\beta
= \modSel_f(N,H)+O_{\varepsilon}(N^{\varepsilon}H^3)\ ,
$$
\par
\noindent 
that is a {\it modified} version of Gallagher's Lemma [Ga, Lemma 1]. 
In order to establish
the aforementioned estimate of $J_3$ the first author in [C5] applies
such a version of Gallagher's Lemma together with the following further property of the essentially bounded and balanced functions; namely (see [C5] for the proof), 
if, for an absolute constant $A\in [0,1)$ and 
for a fixed $\delta\in(0,1/2)$ one has $N^{\delta}\ll H_1\ll H_2\ll N^{1/2-\delta}$ and 
$$
\modSel_f(N,H)\EssBdd NH^{1+A},\enspace \forall H\in [H_1,H_2]\ , 
\enspace \hbox{\rm then} \enspace 
J_f(N,H)\EssBdd NH^{1+{{1+3A}\over {5-A}}},\enspace \forall H\le H_2\ . 
$$
\smallskip
\par
\noindent 
What about the trade of information from the second generation to the first?
\smallskip
\par
Let us take $f=a+b$ real and essentially bounded, with $b$ balanced. Then, we write 
$$
\sum_{h\le H}\Corr_f(h)\buildrel{\hbox{\piccolissimo (T)}}\over{\sim}\sum_{h\le H}\sum_{n\sim N}f(n)f(n+h)
=\avesum f(x)\sum_{x<m\le x+H}f(m)=
$$
$$
=\avesum f(x)\sum_{x<m\le x+H}a(m)+\avesum f(x)\sum_{x<m\le x+H}b(m)\ .
$$
\par
\noindent
Again by the Cauchy inequality one has
$$
\avesum f(x)\sum_{x<m\le x+H}b(m)\ll \Big(\avesum |f(x)|^2\avesum \Big| \sum_{x<n\le x+H}f(n)-M_f(x,H)\Big|^2\Big)^{1/2}
\EssBdd N^{1/2}J_f(x,H)^{1/2}
\ , 
$$
\par
\noindent
where the mean value is 
$$
M_f(x,H)\defineq \sum_{x<n\le x+H}a(n)\ . 
$$
\par
\noindent
Hence, we get exactly an asymptotic formula with main term 
$$
\avesum f(x)M_f(x,H)\ , 
$$
\par
\noindent
whenever the tails are negligible and mostly the remainder term 
$$
\avesum f(x)\sum_{x<m\le x+H}b(m)\EssBdd N^{1/2} J_f(x,H)^{1/2} 
$$
\par
\noindent
turns out to be sufficiently small after halving the exponent's gain on $J_f(x,H)$.
\medskip
\par
Moreover, when $a$ is quasi-constant, one has
$M_f(x,H)\buildrel{\hbox{\piccolissimo (T)}}\over{\sim}a(x)H$, that implies
$$
\sum_{h\le H}\Corr_f(h)\buildrel{\hbox{\piccolissimo (T)}}\over{\sim} 
H\avesum a(x)^2+H\avesum a(x)b(x)+O_{\varepsilon}\big(N^{1/2+\varepsilon}J_f(x,H)^{1/2}\big)\buildrel{\hbox{\piccolissimo (T)}}\over{\sim} 
$$
$$
\buildrel{\hbox{\piccolissimo (T)}}\over{\sim}H\int_N^{2N} a(t)^2{\rm d}t+H\avesum a(x)b(x)
 +O_{\varepsilon}\big(N^{1/2+\varepsilon}J_f(x,H)^{1/2}\big)\ .
$$
\par
\noindent
In particular, if $f$ is stable of exponent $\alpha$, then by applying partial summation 
to $\avesum a(x)b(x)$ one definitively gets the asymptotic  inequality for the {\it deviation} of $f$, say, 
$$
\D_f(N,H)\defineq \sum_{h\le H}\Corr_f(h)-H\avesum a(x)^2 
\EssBdd NHG^{-1/4}+N^{\alpha}H+H^2\ ,
$$
\par
\noindent
whenever a non-trivial estimate for the second generation, $J_f(x,H)\ll NH^2G^{-1/2}$, holds for some ranges of short intervals width,
large enough in terms of the exponent $\alpha$. 
\medskip
\par
Hence, we find a possible general chain of implications of non-trivial estimates as
$$
\modSel_f(N,H)\EssBdd NH^2G^{-1}\Longrightarrow J_f(N,H)\EssBdd NH^2G^{-1/2}\Longrightarrow
\D_f(N,H)\EssBdd NHG^{-1/4}\ .
$$
\par				
\noindent
The exponent's gain halves at each step, but if it remains a neat positive one, 
then we say that $f$ is {\it stable through generations.}
In $\S6$ complete calculations are supplied for the case of the divisor function $d_3$, while
in future papers we are going to explore the hardest case of stable arithmetic functions 
(eventually) having no logarithmic polynomial.

\bigskip

\par
\noindent {\bf 5. Toccata e fuga: the rare cases of asymptotic formul\ae\  for correlations.}
\smallskip
\par
\noindent
An expected asymptotic formula for the single correlations of $f$ usually takes the form 
(compare [BP]) 
$$
\Corr_f(h)=\SingSer_f(h)\SingInt_f(N)+{\cal R}_f(N,h)\ ,
$$
\par
\noindent
where the product of the so-called {\it singular series}
$\SingSer_f$ and the {\it singular integral} $\SingInt_f$ constitutes the main term, while ${\cal R}_f(N,h)$
has to be proved an error term. Such terminology is customary within the Circle Method, 
that was originally introduced by Hardy, Littlewood and Ramanujan in 1918-20 to attack several additive Diophantine problems.
One of the most famous and maddening problem is the infinitude of the pairs of $2h-${\it twin} primes $(n,n-2h)$
that can be formulated in terms of the correlation of the von Mangoldt function:
$$
\Corr_{\Lambda}(2h)=\sum_{n\sim N}\Lambda(n)\Lambda(n-2h)
$$
\par
\noindent
Indeed, with the aid of the powerful analytic tools, 
Hardy and Littlewood predicted that for every sufficiently large $N$ one has
$$
\Corr_{\Lambda}(2h)=
\SingSer_{\Lambda}(2h)N+{\cal R}_{\Lambda}(N,2h)
$$
\par
\noindent
with a certain singular series $\SingSer_{\Lambda}(2h)\gg 1$ and a conjectured remainder term
${\cal R}_{\Lambda}(N,2h)\ll N(\log N)^{-A}$ for every constant $A>0$. 
Note that  such an asymptotic formula would imply the infinitude
of the $2h-$twin primes, but nowadays nobody has yet proved such a conjecture. 
\medskip
\par
\noindent
Apart from very special cases of functions or some trivial instances of the correlations (as they could be when $h=0$), 
the lack of asymptotic formul\ae\ involves the correlations of the most significant arithmetic functions in number theory. 
One of the exceptional cases is given by the well-known {\it binary additive divisor problem}, 
$$
\Corr_{\divisor}(a)=\sum_{n\sim N}\divisor(n)\divisor(n-a) 
$$
\par
\noindent
where $\divisor$ is  the divisor function. This problem was known at least since [E] time and
has been studied mainly through the consolidated theory of modular forms on $SL(\Z,2)$ (see [IwKo], [Vi]). By
adapting Motohashi's results [Mo] to the problem in short intervals, namely $0\not=a=o(N)$, one has
$$
\Corr_{\divisor}(a)=\SingSer_{\divisor}(a)N+{\cal R}_f(N,a)\ , 
$$
\par
\noindent
where
$$
\SingSer_{\divisor}(a)=\SingSer_{\divisor}(a,\log N)\defineq \sum_{i=0}^{2}(\log N)^i\sum_{j=0}^{2}c'_{i,j}\sum_{d|a}{{\log^j d}\over j}\ , 
$$
$$
{\cal R}_{\divisor}(N,a)\EssBdd N^{2/3}+N^{1/2}|a|^{9/40}+|a|^{7/10} 
$$
\par
\noindent
(here these absolute \thinspace $c'_{i,j}$ \thinspace are not Motohashi's constants [Mo, p.$530$]). 
In search of remarkable improvements on the latter formula, one has to be content with an extensive literature on moments of ${\cal R}_{\divisor}(N,a)$
(see [IM]). 
The general case of the additive divisor problem for $d_k$, i.e. establishing an asymptotic formula for 
$$
\Corr_k(a)\defineq \sum_{n\sim N}d_k(n)d_k(n-a)\ , 
$$
\par
\noindent
is much harder and still unsolved when $k\ge 3$.  Even the basic case of $\Corr_3(a)$ seems to be hopeless with present technology due to the poor knowledge about the structure of $SL(\Z,3)$. We address the interested reader to [C2] for a very short tale about this fascinating problem. Moreover, we recall that the Dispersion Method of Linnik (see his book [L]) made its appearance to attack these kind of binary problems.
In $\S2$ we show only an elementary version of the Dispersion Method, while
we have to mention the beautiful paper [Iw] as an example of a non-trivial application of it.
\medskip
\par
\noindent
In the modular forms environment the correlations are widely 
known as {\it Shifted Convolution Sums} (see [Mi]). In such a context,
the possibility of establishing asymptotic formul\ae\ depends directly on the same structure of the modular forms.
This is particularly successful for
the special class of arithmetic functions given by the Hecke eigenvalues $\lambda(n)$. For example, 
Conrey and Iwaniec [CoIw] provide an asymptotic estimate
$$
\Corr_{\lambda}(h)\sim N \sum_{q=1}^{\infty}c_q(h)p(q)^2\ ,
$$
\par
\noindent
where $c_q(h)$ is a Ramanujan sum, while we refer the reader to the quoted paper for the definition of
$p(q)$ and the intricate remainder term. It is remarkable that in a joint and unpublished work with Iwaniec the first author 
has easily deduced from the aforementioned formula of [CoIw] non-trivial bounds for the symmetry integral of the eigenvalues $\lambda$ (essentially square-root cancellation). 
Other spectacular advances on bounding the correlations of the Hecke eigenvalues have been achieved 
by Holowinsky in [Ho1], [Ho2] and then applied jointly with Soundararajan in [HoSo]. 
\medskip
\par
\noindent
Returning to more familiar functions as $d_k$, it seems that the only hope remains the Large Sieve
and all the methods strictly related to such an inequality. 
In fact, similarly to the binary additive divisor problem, the Large Sieve turns out to be crucial when one looks for asymptotic formul\ae\ for
 \lq \lq {\it mixed}\rq \rq correlations of $\divisor$ with some \lq \lq reasonable\rq \rq \thinspace multiplicative arithmetic function (see [C0], where the fundamental reference is Linnik's book [L]). 
Essentially the reason is that, due to the Dirichlet hyperbola method, the divisor function 
has \lq \lq {\it level}\rq \rq \ $\ell<1/2$ of distribution in the arithmetic progressions, which is the same \lq \lq {\it large-sieve-barrier}\rq \rq\ 
for primes in arithmetic progressions, i.e. the celebrated Bombieri-Vinogradov Theorem. Note that the barrier $\ell <1/2$ prevents one from applying
this very strong tool  to the correlations of the von Mangoldt function $\Lambda$, but it is harmless with the \lq \lq truncated von Mangoldt\rq \rq\ function $\Lambda_R$, used in the breakthrough of Goldston-Pintz-Yildirim on small gaps between primes [GPY]. 
An alternate approach  appeals to Duke-Friedlander-Iwaniec [DFI] bounds for bilinear forms of Kloosterman fractions instead of 
the Bombieri-Vinogradov theorem. It allowed the first author [C] to successfully establish asymptotic formul\ae\ for the correlations of essentially bounded 
functions $f$ such that the Eratosthenes transform $f\ast \mu$ is supported up to $O(x^{{1\over 2}+{1\over {190}}-\varepsilon})$ and $\Lambda_R(n)$ with $R\ll x^{{1\over 2}+{1\over {190}}-\varepsilon}$ might be a remarkable example. 
However, this is a 
small improvement on the \lq \lq level\rq \rq\ for correlations, with respect to the aforementioned level in the arithmetic progressions given by the Large Sieve barrier. \enspace See the book [El] for the links between the concepts of {\it level}.
\smallskip
\par
A further possibility is open when $f\ast \mu$ vanishes outside of a very sparse set, where the \lq \lq low density\rq \rq \ has an actual effect 
on the level in arithmetic progressions. In this direction, we refer the reader to [BPW] and [To] 
on the $k-$free numbers, whose characteristic function is defined by 
$$
\sum_{d^k|n}\mu(d)=\sum_{q|n}g_k(q)\ . 
$$
\par
\noindent
Indeed, it is plain that here the Eratosthenes transform $g_k$ is supported on the $k-$th powers, a very \lq \lq low-density\rq \rq \thinspace support. 

\bigskip

\par
\noindent {\bf 6. Crescendo: asymptotic formul\ae\ for $d_3$ in almost all short intervals.}
\smallskip
\par
\noindent
Here we turn our attention to non-trivial bounds for the Selberg integral of $d_3$, postponing the general discussion on $d_k$ to next section. 
Ivi\'c [Iv2] proved that the inequality
$$
J_3(N,H)
\ll N^{1-\varepsilon_1}H^2 
$$
\par
\noindent
holds for the width $\theta>1/6$ with a neat exponent's gain $\varepsilon_1>0$. 
In other words, defining $\theta>\theta_3$ as
the range of the admissible width of the short interval for such an inequality, Ivi\'c has proven $\theta_3=1/6$. 
This is built upon the value $\sigma_3\le 7/12$, where $\sigma_k$ is the so-called {\it Carlson's abscissa} for the 
Riemann zeta $2k$-th moment, 
$$
\sigma_k \defineq \inf \{ \sigma \in [0,1] : \int_{1}^{T}\left| \zeta(\sigma+it)\right|^{2k} {\rm d}t \ll T \}
$$
\par				
\noindent
(see next section for some of the known values $\sigma_k$ quoted from [Iv0]). 
\smallskip
\par
Conjecture CL provides improvements on Ivi\'c's result since for every width $\theta\le 1/3$ it yields 
the \lq \lq best possible estimate\rq \rq, i.e. the square-root cancellation $
\modSel_3(N,H)
\EssBdd NH$,
which in turn implies the optimal $\theta_3=0$ through the arguments of $\S4$. 
This estimate allows improvements on recent results [BBMZ] and [IW] for sums of correlations of $d_3$. Further, it is worth to
recall here that the lower bound $J_3(N,H)\gg NH\log^{4}N$ holds if $0<\theta<1/3$ (see [C0]).  
\smallskip
\par
Now let us prove Corollary 1.
\smallskip
\par
\noindent {\stampatello Proof of Corollary 1}. From  the decomposition $d_3(n)=p_2(\log n)+\widetilde{d_3}(n)$ introduced in $\S3$ one gets
$$
\Corr_3(h)\defineq \sum_{n\sim N}d_3(n)d_3(n-h)
=\sum_{N+h<m\le 2N+h}d_3(m+h)d_3(m)=\sum_{n\sim N}d_3(n+h)d_3(n)+O_{\varepsilon}(N^{\varepsilon}|h|)=
$$
$$
=\sum_{n\sim N}\widetilde{d_3}(n)\widetilde{d_3}(n+h)+2\sum_{n\sim N}\widetilde{d_3}(n)p_2(\log n)+\sum_{n\sim N}p_2(\log n)^2
 + O_{\varepsilon}\left( N^{\varepsilon}|h|\right)\ , 
$$
\par
\noindent
where recall that
$d_3, \widetilde{d_3}$ are essentially bounded, while $p_2(\log n)$ is a quasi-constant function of $n$. 
Therefore,
$$
\D_3(N,H)=\sum_{h\le H}\Corr_3(h)-H\sum_{n\sim N}p_2(\log n)^2=\sum_{n\sim N}\widetilde{d_3}(n)\!\!\sum_{n<m\le n+H}\widetilde{d_3}(m)
+ 2H\sum_{n\sim N}\widetilde{d_3}(n)p_2(\log n)+ O_{\varepsilon}(N^{\varepsilon}H^2). 
$$
\par
\noindent
By applying 
partial summation and $(\widetilde{\ast})$ one has
$$
\sum_{n\sim N}\widetilde{d_3}(n)p_2(\log n)
=p_2(\log (2N))\sum_{n\sim N}\widetilde{d_3}(n)-\int_{N}^{2N}\!\!\sum_{N<n\le t}\widetilde{d_3}(n){{p'_2(\log t)}\over t}{\rm d}t
\EssBdd \max_{t\le 2N}\Big| \sum_{N<n\le t}\widetilde{d_3}(n)\Big|\EssBdd N^{\alpha_3}\ .
$$
\par
\noindent
Since the lower bound $J_3(N,H)\gg NH\log^4N$ holds at least for width $0<\theta<1/3$ (see [C0]), then
Cauchy's inequality implies
$$
\sum_{n\sim N}\widetilde{d_3}(n)\!\!\sum_{n<m\le n+H}\widetilde{d_3}(m)
\EssBdd \sqrt{N\Big(\avesum \Big| \sum_{x<m\le x+H}d_3(m)-M_3(x,H)\Big|^2 + {{H^4L^2}\over N}\Big)}
\EssBdd \sqrt{N J_3(N,H)}\ , 
$$
\par
\noindent
where recall that $L=\log N$ and (see $\S3$) 
$$
M_3(x,H)=Hp_2(\log x)=\sum_{x<m\le x+H}p_2(\log m)+O(x^{-1}H^2 L)\ .
$$
Thus, the conclusion follows immediately from Theorem 1.
\hfill \square

\bigskip

\par
\noindent {\bf 7. Main Theme: from all long intervals to almost all short intervals.}
\smallskip
\par
\noindent
The inductive identity $d_k=d_{k-1}\ast \1$ invites to explore a possible path in order to generalize our conjectures and results to each divisor function $d_k$  for $k>3$ by infering formul\ae\ for $d_k$ in almost all short intervals from suitable information on $d_{k-1}$ in long intervals. However, the actual known values of the amplitudes $\alpha_k$ seem to be a first serious bottle-neck.
Further, while it might be comparatively easy to attack Conjecture CL, the path climbs up drastically when it comes to
the general case $k>3$. Although a general $k-$folding method is available (see next Lemma 2, that is essentially the core of such a method), more and more technical problems are foreseeable as $k$ increases.
Besides, at the outset one has to face the problem of showing sufficient proximity of the analytic and the arithmetic forms of the mean value, at least by mean-square approximation, i.e. an inequality of the form
$$
\avesum \left| M_k(x,H)-\widetilde{M}_k(x,H)\right|^2 \EssBdd N^{1-2/k}H^2\ ,
$$
\par				
\noindent
for a suitable choice of $\widetilde{M}_k(x,H)$. Being unconceivable to give a rigorous general proof of such a proximity by direct calculations as we did for the case $k=3$ in \S 3, at the moment  we have to content ourself with the following heuristic considerations.
Some results of Ivi\'c [Iv2] provide non-trivial estimates of the Selberg integral $J_k(N,H)$ with the mean value $M_k(x,H)$ assigned in the appropriate analytic form.
Hence, Theorem 2 legitimates the assumption that the  arithmetic mean value $\widetilde M_k(x,H)$
is close to the analytic counterpart $M_k(x,H)$ for every $k>3$
in the ranges of the short interval width $\theta$ provided by Ivi\'c's results.
Actually, we know that the analytic form is $H$ times a $k-1$ degree polynomial in $\log x$, the same shape that approximates (see next Lemma 2) the arithmetic form. Then, comparing Ivi\'c's results with ours in a common range for $\theta$, we easily conclude that these polynomials must coincide 
(Amitsur formula [A] gives polynomials' degree and Tull's Lemma [Tu] the remainders). 
\bigskip
\par
\noindent
Now, let us turn our attention to next Lemma 2, where we show the so-called \lq \lq $k-$folding method\rq \rq. At this aim,
assuming that $H$ does not depend on $x$ explicitly, 
let us consider for any fixed $k\ge 2$ the weighted sum 
$$
S_k(x,H)\defineq \sum_n a_k(n)w(n-x) 
= \multiplesum_{n_1 \enspace \qquad \enspace n_k}a(n_1)\cdots a(n_k)w(n_1\cdots n_k-x)\ ,
$$
\par
\noindent
where $a_k\defineq \underbrace{a\ast \cdots \ast a}_{k\, {\rm times}}$ is the $k-$fold Dirichlet product of an arithmetic function $a:\N \rightarrow \C$ and the weight $w:\N \rightarrow \C$ is supposed to be uniformly bounded in its support, that is contained in $[-H,H]$.
\par
\noindent
Further, for $M\defineq [(N-H)^{1/k}]$ let us denote
$$
g_k(q)=g_k(q,a,M)\defineq 
a_{k-1}(q)+\sum_{j\le k-1}a_{k-1}^{(j)}(q)\ , 
$$
\par
\noindent
where
$$
a_{k-1}^{(j)}(q)\defineq \multiplesum_{{n_1 \enspace \qquad \enspace n_{k-1}}\atop {{n_1 \cdots n_{k-1} = q}\atop {n_1,\ldots,n_j<M}}}a(n_1)\cdots a(n_{k-1})\qquad \forall j\le k-1\ .
$$
\medskip
\par
\noindent {\stampatello Lemma 2.} {\it If the arithmetic function $a$ is essentially bounded, then $\forall \varepsilon>0$ and uniformly for every} $x\sim N$, 
$$
S_k(x,H)=\sum_{q\le {x/M}}g_k(q)\sum_{{0\le |n-x|\le H}\atop {n\equiv 0(\!\!\bmod q)}}a\Big({n\over q}\Big)w(n-x)
+ O_{k,\varepsilon}\Big( N^{\varepsilon}\Big({H\over {N^{1/k}}}+{{H^2}\over N}+1\Big)\Big)\ . 
$$
\smallskip
\par
\noindent {\stampatello Proof.}$\!$ First note that $n_1 \cdots n_k\ge x-H\ge N-H$ implies that
at least one of $n_1, \ldots , n_k$ has to be $\ge M$. Then, let us define the partial sums $\Sigma_0,\Sigma_1,\ldots,\Sigma_{k-1}$ of $S_k(x,H)$ 
as follows:
\par
$\Sigma_0$ is the part of $S_k(x,H)$ corresponding to $n_1\ge M$\ , 
\par
$\Sigma_1$ is the part of $S_k(x,H)-\Sigma_0$ corresponding to $n_2\ge M$, 
\par
$\Sigma_2$ is the part of $S_k(x,H)-\Sigma_0-\Sigma_1$ corresponding to $n_3\ge M$,
and so on.
\par
\noindent 
Therefore, we set $a_{k-1}^{(0)}(q)\defineq a_{k-1}(q)$ and split $S_k(x,H)$ as
$$
S_k(x,H)=\Sigma_0+\Sigma_1+\ldots+\Sigma_{k-1}=
$$
$$
=\sum_{q\le {{x+H}\over M}}a_{k-1}^{(0)}(q)\sum_{{n=qm}\atop {m\ge M}}a(m)w(n-x) + \sum_{q\le {{x+H}\over M}}a_{k-1}^{(1)}(q)\sum_{{n=qm}\atop {m\ge M}}a(m)w(n-x) + 
$$
$$
+ \cdots + \sum_{q\le {{x+H}\over M}}a_{k-1}^{(k-1)}(q)\sum_{{n=qm}\atop {m\ge M}}a(m)w(n-x)=\sum_{q\le {{x+H}\over M}}g_k(q)\sum_{{n=qm}\atop {m\ge M}}a(m)w(n-x) \ ,
$$
\par
\noindent
where for each $j=0,1,\ldots,k-1$ the multiple sum
$$
a_{k-1}^{(j)}(q)=\multiplesum_{{n_1 \enspace \qquad \enspace n_{k-1}}\atop {{n_1 \cdots n_{k-1} = q}\atop {n_1,\ldots,n_j<M}}}a(n_1)\cdots a(n_{k-1})
$$
\par				
\noindent
has $j$ variables restricted by $M$ (that depends on $N,H,k$, but not on $x$).
\par
\noindent
Observe that, since
$\displaystyle{g_k(q)=\sum_{j=0}^{k-1}a_{k-1}^{(j)}(q)}$, then 
$|g_k(q)|\le k{|a|}_{k-1}(q)\EssBdd 1$, where we set 
${|a|}_{k-1}\defineq \underbrace{|a|\ast \cdots \ast |a|}_{k-1\, {\rm times}}$. Thus, 
$$
S_k(x,H)
= \sum_{q\le {{x+H}\over M}}g_k(q)\sum_{{{{x-H}\over q}\le m\le {{x+H}\over q}}\atop {m\ge M}}a(m)w(qm-x) = 
$$
$$
= \sum_{q\le {{x-H}\over M}}g_k(q)\sum_{{{x-H}\over q}\le m\le {{x+H}\over q}}a(m)w(qm-x) 
+ {\cal O}_{k,\varepsilon}\Big( \sum_{{{x-H}\over M}<q\le {{x+H}\over M}}\sum_{\left|m-{x\over q}\right|\le {H\over q}}x^{\varepsilon}\Big) = 
$$
$$
= \sum_{q\le {x/M}}g_k(q)\sum_{{{x-H}\over q}\le m\le {{x+H}\over q}}a(m)w(qm-x) 
+ {\cal O}_{k,\varepsilon}\Big( \sum_{{{x-H}\over M}<q\le {{x+H}\over M}}\sum_{\left|m-{x\over q}\right|\le {H\over q}}x^{\varepsilon}\Big) = 
$$
$$
= \sum_{q\le {x/M}}g_k(q)\sum_{{0\le |n-x|\le H}\atop {n\equiv 0(\!\!\bmod q)}}a\left({n\over q}\right)w(n-x) 
+ {\cal O}_{k,\varepsilon}\Big( N^{\varepsilon}\sum_{{{x-H}\over M}<q\le {{x+H}\over M}}\sum_{\left|m-{x\over q}\right|\le {H\over q}}1\Big)\ .
$$
\par
\noindent
Since \thinspace $\displaystyle{q>{{x-H}\over M}\gg {N\over M}}$ \thinspace uniformly as $N\le x\le 2N$, then the ${\cal O}_{k,\varepsilon}$-term is 
$$
\ll N^{\varepsilon}\sum_{{{x-H}\over M}<q\le {{x+H}\over M}}\left( {H\over q} + 1\right) 
\ll N^{\varepsilon}\left( {H\over M} + 1\right) \left( {{HM}\over N} + 1\right) 
\ll N^{\varepsilon}\Big({{H^2}\over N} + {H\over {N^{1/k}}} + 1\Big)\ .\enspace \square
$$
\medskip
\par
\noindent
An application of Lemma 2 to attack the General Conjecture CL would require $a=\1$, so that $a_k=d_k$, attached to
the Cesaro weights.
\medskip
\par
\noindent {\bf Remark.} The General Conjecture CL together with Theorem 2 justifies the title of the present section:  
we positively interpret the (upper) bounds for the amplitude $\alpha_{k-1}$ as a property which holds in 
every long interval, while the consequent bound of $\modSel_k$ clearly means an almost all short interval property.
However, since the present state of knowledge assigns values to $\alpha_k$ 
that grow towards $1$ with $k$, the quality of the General Conjecture CL  and Theorem 2 ruins unrelentingly for large values of $k$. 
Of course, as already seen when $k=3$, the general scheme of implications of non-trivial estimates discussed in $\S4$ still works in suitable ranges of short intervals, namely when the width of the short interval is sufficiently large with respect to $\alpha_k$. In particular, one easily proves Corollary 2 on the deviation of $d_k$ (see the definition in \S0) for any $k>3$ by closely following the Corollary 1 proof with the aid of the arguments in $\S4$.
\medskip
\par
\noindent
Note that we have an improvement on the aforementioned Ivi\'c's results on the Selberg integral $J_k(N,H)$ for $k>3$. More precisely,
Theorem 1 in [Iv2] holds for $\theta\in(\theta_k,1)$ with $\theta_k\defineq 2\sigma_k-1$. In particular, 
from Corollary 2 of [Iv2] one has
$$
\theta_4={1\over 4},\qquad \theta_5={{11}\over {30}},\qquad \theta_6={3\over 7}\ , 
$$
\par
\noindent
respectively corresponding to the following known upper bounds [Iv0] for Carlson's abscissae, defined in $\S6$,
$$
\sigma_4 \le {5\over 8},\qquad \sigma_5 \le {{41}\over {60}},\qquad \sigma_6 \le {5\over 7}\ . 
$$
\par
\noindent
Due to Corollary 2 such values of $\theta_k$ 
are superseded by  
$$
\widetilde{\theta}_4={{11}\over {128}},\qquad \widetilde{\theta}_5={1\over 5},\qquad \widetilde{\theta}_6={1\over 3}\ , 
$$
\par
\noindent
which follow from the well-known upper bounds for the amplitudes (see [Ti])
$$
\alpha_3 \le {{43}\over {96}}\ (\hbox{\rm Kolesnik}, 1981),
\quad 
\alpha_4 \le {1\over 2}\ (\hbox{\rm Hardy-Littlewood [HL]} ,1922),
\quad
\alpha_5 \le {{11}\over {20}}\ (\hbox{\rm [Iv0], Ch.13})\ .
$$
\par				
\noindent
It might be possible that $\theta_k<\widetilde{\theta}_k$ for any $k>6$, although we think that it is unlikely.  
However, we put our
emphasis essentially on the method: unlike Ivi\'c, we
do not use any deep property of 
the Riemann $\zeta-$function (at least not explicitly)
and rely uniquely upon known values 
of the exponents $\alpha_k$; but one should not forget that most known values $\alpha_k$ follow directly from non-trivial estimates of the moments of the Riemann zeta function on the line (for a clear digression on such a topic see the wonderful book by Ivi\'c [Iv0]). On the other side, as the first author pointed out (see [C2], Theorem 1.1), estimates of the Selberg integral of $d_k$ have non-trivial consequences  on the $2k-$th moments of the Riemann zeta function.
Thus, recalling that Kolesnik has found his bound for $\alpha_3$ without the aid of the $6-$th moments and that $\alpha_k$ known values 
go to $1$, it is plain that, at least for relatively low values of $k$, Theorem 2 confirms that 
the Selberg integrals of $d_k$ and the $2k-$th moments of the Riemann zeta function 
are connected by a circle route. Next section is devoted to a further discussion on the argument.

\bigskip

\par
\noindent {\bf 8. Finale: conditional bounds for the moments of $\zeta$ on the critical line.}
\smallskip
\par
\noindent
As already mentioned in the previous section, at least in theory we could draw some consequences of Theorem 2
on the moments of the Riemann zeta function. However, our method does not lead to
a better result than those available in the literature for the 
$2k-$th moment when $4\le k\le 6$ and maybe the scenario is even worst  when $k>6$.
The reason is essentially the bound $N^{1-1/k}H^2$ in the General Cojecture CL, that is an unavoidable barrier term. It transfers
from $\modSel_k(N,H)$ to a bound of the $2k-$th moment,
via the Selberg integral $J_k(N,H)$, whenever we appeal to Theorem 1.1 of [C2].
\par
Here we take the opportunity of applying Theorem 1.1 of [C2] to give some conditional bounds,
depending on estimates for the Selberg integrals $J_k$ which are proved for $k=3$, but unproved for $k>3$. 
At this aim, let us define the {\it excess} $E_k$ as a real number such that
$$
I_k(T)=\int_{T}^{2T}\Big| \zeta\Big( {1\over 2}+it\Big)\Big|^{2k} {\rm d}t \EssBdd T^{1+E_k}\ .
$$
\par
\noindent
With an application of H\"older's inequality
from the well-known values $E_2=0$ and $E_6=1$ (from Heath-Brown's bound $I_6(T)\EssBdd T^2$) one gets 
$$
E_3=1/4, \quad E_4=1/2, \quad E_5=3/4\ .
$$
\par
\noindent
According to this definition, Theorem 3 gives $E_k={k\over 2}(A+B)-B$ whenever $\displaystyle{J_k(N,H)\EssBdd N^{1+A}H^{1+B}}$ holds
for $H\ll N^{1-2/k}$. In particular, as already mentioned, it gives $E_3=1/10$, when combined with recent [C5]. 
\medskip
\par
\noindent {\stampatello Proof of Theorem 3.} In [C2] the Selberg integral of $d_k$ is defined in the standard way as
$$
\int_{hx^{\varepsilon}}^{x}\Big|\sum_{t<n\le t+h}d_k(n)- M_k(t,h)\Big|^2 {\rm d}t\ . 
$$
\par
\noindent
It is easy to see that a dyadic argument allows to replace 
$[hx^{\varepsilon},x]$ by dyadic intervals like $[N,2N]$, and the substitution of 
the integral on $[N,2N]$ by $J_k(N,H)$ generates only negligible remainder terms. Hence, 
we may apply Theorem 1.1 of [C2] by using $J_k(N,H)$ instead of the above integral and $\forall \varepsilon>0$ small we get
$$
I_k(T)\EssBdd T+T\max_{T^{1+\varepsilon}\ll N\ll T^{k\over 2}}{T\over {N^2}}\max_{0<H\ll {N\over T}}J_k(N,H)
\EssBdd T+T\max_{T^{1+\varepsilon}\ll N\ll T^{k\over 2}}N^A\left({N\over T}\right)^B\EssBdd 
$$
$$
\EssBdd T+T\left(T^{k\over 2}\right)^A\left(T^{{k\over 2}-1}\right)^B
\EssBdd T^{1+{k\over 2}(A+B)-B}\ .\enspace \square
$$

\bigskip

\par
\noindent {\bf 9. Epilogo: the best unconditional exponent of $|\zeta|^6$ mean on the line.}
\smallskip
\par
\noindent
As we saw in the introduction, a further immediate consequence of the recent result [C5] is that
we have an elementary deduction of $E_3=1/10$ just taking
$A=0$ and $B=1/5$ in Theorem 3 for $H\ll N^{1/3}$. 
However, the proof of Theorem 1.1 [C2] uses the approximate functional equation for $\zeta^k$ (and 
Gallagher's Lemma), 
so the whole study actually proving $E_3=1/10$ (in [C5]
Gallagher's Lemma is also applied) is just shorter, rather than elementary, with respect to the Heath-Brown's bound of $I_6$. 
\par
Furthermore, see Corollary of [C5], our new approach based on the \lq \lq modified Gallagher Lemma\rq \rq, contained in the 
forthcoming paper [CL], assures under Conjecture CL the even better excess $E_3=0$ for the Riemann zeta function, namely 
the well-known {\it weak $6-$th moment}. 
\par
\noindent
The excesses $E_k$ here ($k\ge 3$) are from a $2$nd generation approach, while Ivi\'c's [Iv1] is a $1$st generation one. 

\bigskip

\par
\noindent
\centerline{\stampatello References}
\smallskip
\item{\bf [A]} Amitsur, S. A. \thinspace - \thinspace {\sl Some results on arithmetic functions} \thinspace - \thinspace  J. Math. Soc. Japan {\bf 11} (1959), 275--290. $\underline{\tt MR\thinspace 26 \#67}$ \thinspace - \thinspace available online 
\item{\bf [BBMZ]} Baier, S., Browning, T.D., Marasingha, G. and Zhao, L. \thinspace - \thinspace {\sl Averages of shifted convolutions of $d_3(n)$} \thinspace - \thinspace http://arxiv.org/abs/1101.5464v2
\item{\bf [BP]} Br\"udern, J., Perelli, A. \thinspace - \thinspace {\sl A note on the distribution of sumsets} \thinspace - \thinspace Funct. Approx. Comment. Math. {\bf 29} (2001), 81--88. $\underline{\tt MR\thinspace 2005m\!:\!11033}$ 
\item{\bf [BPW]} Br\"udern, J., Perelli, A. and Wooley, T. \thinspace - \thinspace {\sl Twins of k-Free Numbers and Their Exponential Sum} \thinspace - \thinspace Michigan Math. J. {\bf 47} (2000), No. {\bf 1}, 173--190. $\underline{\tt MR\thinspace 2001c\!:\!11113}$ \thinspace - \thinspace available online 
\item{\bf [CoIw]} Conrey, B. and Iwaniec, H. \thinspace - \thinspace {\sl Spacing of zeros of Hecke $L$-functions and the class number problem} \thinspace - \thinspace Acta Arith. {\bf 103} (2002), no. {\bf 3}, 259--312. $\underline{\tt MR\thinspace 2003h\!:\!11103}$ 
\item{\bf [C]} Coppola, G. \thinspace - \thinspace {\sl On the Correlations, Selberg integral and symmetry of sieve functions in short intervals, III} \thinspace - \thinspace http://arxiv.org/abs/1003.0302v1 
\item{\bf [C0]} Coppola, G. \thinspace - \thinspace {\sl On some lower bounds of some symmetry integrals} \thinspace - \thinspace http://arxiv.org/abs/1003.4553v2 - to appear on Afrika Mathematika (Springer)  
\item{\bf [C1]} Coppola, G. \thinspace - \thinspace {\sl On the modified Selberg integral} \thinspace - \thinspace http://arxiv.org/abs/1006.1229v1 
\item{\bf [C2]} Coppola, G. \thinspace - \thinspace {\sl On the Selberg integral of the $k$-divisor function and the $2k$-th moment of the Riemann zeta-function} \thinspace - \thinspace Publ. Inst. Math. (Beograd) (N.S.) {\bf 88(102)} (2010), 99--110. \thinspace - \thinspace available online 
\item{\bf [C3]} Coppola, G.\thinspace - \thinspace {\sl On the Correlations, Selberg integral and symmetry of sieve functions in short intervals, II} \thinspace - \thinspace Int. J. Pure Appl. Math. {\bf 58.3}(2010), 281--298. \thinspace - \thinspace available online 
\item{\bf [C4]} Coppola, G. \thinspace - \thinspace {\sl On the symmetry of divisor sums functions in almost all short intervals} \thinspace - \thinspace Integers {\bf 4} (2004), A2, 9 pp. (electronic). $\underline{\tt MR\enspace 2005b\!:\!11153}$ 
\item{\bf [C5]} Coppola, G. \thinspace - \thinspace {\sl On the Selberg integral of the three-divisor function $d_3$} \thinspace - \thinspace available online at the address http://arxiv.org/abs/1207.0902v3 
\item{\bf [CL]} Coppola, G. and Laporta, M. \thinspace - \thinspace {\sl A modified Gallagher's Lemma} \thinspace - \thinspace http://arxiv.org/abs/ 
\item{\bf [CS]} Coppola, G. and Salerno, S. \thinspace - \thinspace {\sl On the symmetry of the divisor function in almost all short intervals} \thinspace - \thinspace Acta Arith. {\bf 113} (2004), no. {\bf 2}, 189--201. $\underline{\tt MR\thinspace 2005a\!:\!11144}$ 
\item{\bf [D]} Davenport, H. \thinspace - \thinspace {\sl Multiplicative Number Theory} \thinspace - \thinspace Third Edition, GTM 74, Springer, New York, 2000. $\underline{{\tt MR\thinspace 2001f\!:\!11001}}$ 
\item{\bf [De]} De Roton, A. \thinspace - \thinspace {\sl On the mean square of the error term for an extended Selberg class} \thinspace - \thinspace Acta Arith. {\bf 126} (2007), no. {\bf 1}, 27--55. $\underline{\tt MR\thinspace 2007j\!:\!11121}$ 
\item{\bf [DFI]} Duke, W., Friedlander, J. and Iwaniec, H. \thinspace - \thinspace {\sl Bilinear forms with Kloosterman fractions} \thinspace - \thinspace Invent. Math. {\bf 128} (1997), no. {\bf 1}, 23--43. $\underline{\tt MR\thinspace 97m\!:\!11109}$ 
\item{\bf [E]} Estermann, T. \thinspace - \thinspace {\sl \"Uber die Darstellungen einer Zahl als Differenz von zwei Produkten} \thinspace - \thinspace J. Reine Angew. Math. (Crelle Journal) {\bf 164} (1931), 173--182. 
\item{\bf [El]} Elliott, P.D.T.A. \thinspace - \thinspace {\sl On the correlation of multiplicative and the sum of additive arithmetic functions} Mem. Amer. Math. Soc. {\bf 112} (1994), no. 538, viii+88 pp. $\underline{\tt MR\thinspace 95d\!:\!11099}$ 
\item{\bf [F]} Fej\'er, L. \thinspace - \thinspace {\sl \"Uber trigonometrische Polynome} \thinspace - \thinspace J. f\"ur Math. (Crelle Journal) {\bf 146} (1916), 53--82. 
\item{\bf [Ga]} Gallagher, P. X. \thinspace - \thinspace {\sl A large sieve density estimate near $\sigma =1$} \thinspace - \thinspace Invent. Math. {\bf 11} (1970), 329--339. $\underline{\tt MR\enspace 43\# 4775}$ 
\item{\bf [Go]} Goldmakher, L. \thinspace - \thinspace {\sl Character sums to smooth moduli are small} \thinspace - \thinspace Canad. J. Math. {\bf 62} (2010), no. {\bf 5}, 1099--1115. $\underline{{\tt MR\thinspace 2011k\!:\!11108}}$ 
\item{\bf [Gr]} Green, B. \thinspace - \thinspace {\sl On arithmetic structures in dense sets of integers} \thinspace - \thinspace Duke Math. J. {\bf 114} (2002), no. {\bf 2}, 215--238. $\underline{{\tt MR\thinspace 2003i\!:\!11021}}$ 
\item{\bf [GPY]} Goldston, D.A., Pintz, J. and Yildirim, C. \thinspace - \thinspace {\sl Primes in tuples. I} \thinspace - \thinspace Ann. of Math. (2) {\bf 170} (2009), no. {\bf 2}, 819--862. $\underline{{\tt MR\thinspace 2011c\!:\!11146}}$ 
\item{\bf [HL]} Hardy, G.H. and Littlewood, J.E. \thinspace - \thinspace {\sl The approximate functional equation in the theory of the zeta-function, with applications to the divisor problems of Dirichlet and Piltz} \thinspace - \thinspace Proc. London Math. Soc.(2) {\bf 21} (1922), 39--74. 
\item{\bf [Ho1]} Holowinsky, R. \thinspace - \thinspace {\sl A sieve method for shifted convolution sums} \thinspace - \thinspace Duke Math. J. {\bf 146}  (2009), no. {\bf 3}, 401--448. $\underline{\tt MR\thinspace 2010b\!:\!11127}$ 
\item{\bf [Ho2]} Holowinsky, R. \thinspace - \thinspace {\sl Sieving for mass equidistribution} \thinspace - \thinspace Ann. of Math. (2) {\bf 172}  (2010), no. {\bf 2}, 1499--1516. $\underline{\tt MR\thinspace 2011i\!:\!11060}$ 
\item{\bf [HoSo]} Holowinsky, R. and Soundararajan, K. \thinspace - \thinspace {\sl Mass equidistribution for Hecke eigenforms} \thinspace - \thinspace Ann. of Math. (2) {\bf 172} (2010), no. {\bf 2}, 1517--1528. $\underline{\tt MR\thinspace 2011i\!:\!11061}$ 
\item{\bf [Iv0]} Ivi\'c, A. \thinspace - \thinspace {\sl The Riemann Zeta Function} \thinspace - \thinspace John Wiley \& Sons, New York, 1985. (2nd ed., Dover, Mineola, N.Y. 2003). $\underline{\tt MR\thinspace 87d\!:\!11062}$ 
\item{\bf [Iv1]} Ivi\'c, A. \thinspace - \thinspace {\sl The general additive divisor problem and moments of the zeta-function} \thinspace - \thinspace New trends in probability and statistics, Vol. {\bf 4} (Palanga, 1996), 69--89, VSP, Utrecht, 1997. $\underline{\tt MR\thinspace 99i\!:\!11089}$ 
\item{\bf [Iv2]} Ivi\'c, A. \thinspace - \thinspace {\sl On the mean square of the divisor function in short intervals} \thinspace - \thinspace J. Th\'eor. Nombres Bordeaux {\bf 21} (2009), no. {\bf 2}, 251--261. $\underline{\tt MR\thinspace 2010k\!:\!11151}$ 
\item{\bf [IM]} Ivi\'c, A. and Motohashi, Y. \thinspace - \thinspace {\sl On some estimates involving the binary additive divisor problem}, Quart. J. Math. Oxford Ser. (2) {\bf 46}, no. {\bf 184} (1995), 471--483. $\underline{\tt MR\thinspace 96k\!:\!11117}$ 
\item{\bf [IW]} Ivi\'c, A. and Wu, J. \thinspace - \thinspace {\sl On the general additive divisor problem} \thinspace - \thinspace http://arxiv.org/abs/1106.4744v2
\item{\bf [Iw]} Iwaniec, H. \thinspace - \thinspace {\sl Almost-primes represented by quadratic polynomials} \thinspace - \thinspace Invent. Math. {\bf 47} (1978), no. {\bf 2}, 171--188. $\underline{\tt MR\thinspace 58 \#5553}$ 
\item{\bf [IwKo]} Iwaniec, H. and Kowalski, E. \thinspace - \thinspace {\sl Analytic Number Theory} \thinspace - \thinspace American Mathematical Society Colloquium Publications, 53. AMS, Providence, RI, 2004. $\underline{\tt MR\thinspace 2005h\!:\!11005}$ 
\item{\bf [KP]} Kaczorowski, J. and Perelli, A. \thinspace - \thinspace {\sl On the distribution of primes in short intervals} \thinspace - \thinspace J. Math. Soc. Japan {\bf 45}  (1993), no. {\bf 3}, 447--458. $\underline{\tt MR\thinspace 94e\!:\!11100}$ 
\item{\bf [KP(012)]} Kaczorowski, J. and Perelli, A. \thinspace - \thinspace {\sl On the structure of the Selberg class, VII: $1 < d < 2$} \thinspace - \thinspace Ann. of Math. (2) {\bf 173} (2011), no. {\bf 3}, 1397--1441.
\item{\bf [L]} \thinspace Linnik, Ju.V.\thinspace - \thinspace {\sl The Dispersion Method in Binary Additive Problems} \thinspace - \thinspace Translated by S. Schuur \thinspace - \thinspace American Mathematical Society, Providence, R.I. 1963. $\underline{\tt MR\enspace 29\# 5804}$
\item{\bf [Mi]} Michel, P. \thinspace - \thinspace {\sl On the Shifted Convolution Problem} \thinspace - \thinspace available online at the following web address, http://tan.epfl.ch/files/content/sites/tan/files/PhMICHELfiles/Fields2003.pdf 
\item{\bf [Mo]} Motohashi, Y. \thinspace - \thinspace {\sl The binary additive divisor problem} \thinspace - \thinspace Ann. Sci. \'Ecole Norm. Sup. (4) {\bf 27} (1994), no. {\bf 5}, 529--572. $\underline{\tt MR\thinspace 95i\!:\!11104}$ 
\item{\bf [S]} Selberg, A. \thinspace - \thinspace {\sl On the normal density of primes in small intervals, and the difference between consecutive primes} \thinspace - \thinspace Arch. Math. Naturvid. {\bf 47} (1943), no. {\bf 6}, 87--105. $\underline{\tt MR\thinspace 7,48e}$ 
\item{\bf [Te]} Tenenbaum, G. \thinspace - \thinspace {\sl Introduction to Analytic and Probabilistic Number Theory} \thinspace - \thinspace Cambridge Studies in Advanced Mathematics, {\bf 46}, Cambridge University Press, 1995. $\underline{\tt MR\thinspace 97e\!:\!11005b}$ 
\item{\bf [Ti]} Titchmarsh, E. C. \thinspace - \thinspace {\sl The theory of the Riemann zeta-function} \thinspace - \thinspace Second edition. Edited and with a preface by D. R. Heath-Brown. {\sl The Clarendon Press, Oxford University Press, New York}, 1986. \hfil $\underline{\tt MR\thinspace 88c\!:\!11049}$ 
\item{\bf [To]} Tolev, D. \thinspace - \thinspace {\sl On the distribution of $r$-tuples of square-free numbers in short intervals} \thinspace - \thinspace Int. J. Number Theory {\bf 2} (2006), no. {\bf 2}, 225–-234. $\underline{\tt MR\thinspace 2008a\!:\!11111}$ 
\item{\bf [Tu]} Tull, J. P. \thinspace - \thinspace {\sl Average order of arithmetic functions} \thinspace - \thinspace Illinois J. Math. {\bf 5} (1961), 175--181. $\underline{\tt MR\thinspace 22 \#10943}$
\item{\bf [Vi]} Vinogradov,  A. I. \thinspace - \thinspace {\sl A generalized square of the zeta function. Spectral decomposition} (Russian)  \thinspace - \thinspace Zap. Nauchn. Sem. S.-Peterburg. Otdel. Mat. Inst. Steklov. (POMI) {\bf 322} (2005), Trudy po Teorii Chisel, 17--44, 251; translation in J. Math. Sci. (N. Y.) {\bf 137} (2006), no. {\bf 2}, 4617--4633. $\underline{{\tt MR\thinspace 2006g\!:\!11175}}$ 
\bigskip
\par

\leftline{\tt Giovanni Coppola\spaziolungo \spaziolungo \qquad \qquad \enspace \thinspace Maurizio Laporta}
\leftline{\tt Universit\`a degli Studi di Salerno\spaziolungo \thinspace Universit\`a degli Studi di Napoli}
\leftline{\tt Home address \negthinspace : \negthinspace Via Partenio \negthinspace 12 \negthinspace -\spaziolungo Dipartimento di Matematica e Appl.}
\leftline{\tt - 83100, Avellino(AV), ITALY\spaziolungo \qquad \qquad \qquad \qquad \enspace \thinspace Compl.Monte S.Angelo}
\leftline{\tt e-page : $\! \! \! \! \! \!$ www.giovannicoppola.name\qquad \qquad \qquad \qquad \qquad \quad \thinspace Via Cinthia - 80126, Napoli, ITALY}
\leftline{\tt e-mail : gcoppola@diima.unisa.it\spaziolungo \qquad \enspace \thinspace e-mail : mlaporta@unina.it}

\bye